\documentclass{article}
\usepackage[utf8]{inputenc}
\usepackage[T1]{fontenc}
\usepackage[english]{babel}
\usepackage{amsthm}
\usepackage{arxiv}
\usepackage{algorithm2e}
\usepackage{authblk}
\newtheorem{theorem}{Theorem}[section]

\newtheorem{definition}{Definition}[subsection]

\usepackage{amsmath,bm}
\usepackage{amsfonts}
\usepackage{dsfont}
\usepackage{nicefrac}
\usepackage{amssymb}
\usepackage{color}
\usepackage{cases}
\usepackage{pifont}
\usepackage{array}
\usepackage{ulem}
\usepackage{empheq}
\usepackage{mathtools}
\usepackage{appendix}
\usepackage{soul}
\usepackage{tcolorbox}
\usepackage{empheq}
\usepackage{lineno}
\usepackage{musicography}
\usepackage{accents}

\DeclarePairedDelimiterX{\inp}[2]{\langle}{\rangle}{#1, #2}

\DeclareMathOperator*{\argmin}{arg\,min}

\newcommand{\RNum}[1]{\uppercase\expandafter{\romannumeral #1\relax}}

\newcommand{\bfx}{{\bf x}}
\newcommand{\bfu}{{\bf u}}
\newcommand{\bfv}{{\bf v}}
\newcommand{\bfw}{{\bf w}}
\newcommand{\bfPhi}{{\bf \Phi}}
\newcommand{\bfh}{{\bf h}}
\newcommand{\bfeta}{{\bf \eta}}
\newcommand{\bfxi}{{\bf \xi}}
\newcommand{\bfX}{{\bf X}}
\newcommand{\bfe}{{\bf e}}
\newcommand{\bfp}{{\bf p}}
\newcommand{\bfq}{{\bf q}}
\newcommand{\bfy}{{\bf y}}

\newcommand{\uflat}{{\bf u}^{\musFlat{}}}
\newcommand{\vflat}{{\bf v}^{\musFlat{}}}
\newcommand{\wflat}{{\bf w}^{\musFlat{}}}

\newcommand{\hflat}{{\bf h}^{\musFlat{}}}
\newcommand{\etaflat}{{\bf \eta}^{\musFlat{}}}
\newcommand{\xiflat}{{\bf \xi}^{\musFlat{}}}

\theoremstyle{definition}
\newtheorem{example}{Example}[subsection]
\usepackage{ulem}
\RestyleAlgo{ruled}

\title{Alignment of Geophysical Fields: a differential geometry perspective}
\author[1]{Yicun Zhen\thanks{Corresponding author: zhenyicun@proton.me}}
\author[2,3]{Valentin Resseguier}
\author[4]{Bertrand Chapron}
\affil[1]{Department of Oceanography, Hohai University, Nanjing, Jiangsu, China}
\affil[2]{LAB, SCALIAN DS, Rennes, France}
\affil[3]{INRAE, OPAALE, Rennes, France}
\affil[4]{Laboratoire d’Océanographie Physique et Spatiale, Ifremer, Plouzan\'{e}, France}

\begin{document}

\maketitle

\begin{abstract}

To estimate the displacements of physical state variables, the physics principles that govern the state variables must be considered. Technically, for a certain class of state variables, each state variable is associated to a tensor field. Ways displacement maps act on different state variables will then differ according to their associated different tensor field definitions. Displacement procedures can then explicitly ensure the conservation of certain physical quantities (total mass, total vorticity, total kinetic energy, etc.), and a  differential-geometry-based optimisation formulated. Morphing with the correct physics, it is reasonable to apply the estimated displacement map to unobserved state variables, as long as the displacement maps are strongly correlated. This leads to a new nudging strategy using all-available observations to infer displacements of both observed and unobserved state variables. Using the proposed nudging method before applying ensemble data assimilation, numerical results show improved preservation of the intrinsic structure of underlying physical processes. 

\end{abstract}

\section{introduction}  %% \introduction[modified heading if necessary]
Since chaotic divergence is an intrinsic property of turbulent geophysical systems, data assimilation is often necessary to improve model trajectories of state variables. Numerous strategies have then been proposed, including nudging methods, 3-Dimensional variational method \cite{Lorenc1986AnalysisMF}, Kalman filter based methods \cite{Houtekamer2001ASE,Houtekamer1998DataTechnique,Ott2004AAssimilation,Kalman1960}, 4-Dimensional variational methods \cite{Courtier4DVar1994}, and particle filters \cite{ParticleFilterDelMoral1997}. Suppose that the background state estimate $x^b(t,\bfx)$ is governed by equation 
\begin{align}
\frac{\partial x^b(t,\bfx)}{\partial t} = M(t, x^b(t,\bfx)),
\end{align}
where $M$ is the numerical model and $\bfx\in\Omega$ is the coordinate of a point in the domain $\Omega$.  The nudging method in theory assumes that the observations $y^o(t)$ are taken in continuous time. The method aims to find matrix $K$, which could depend on time $t$,  background estimate $x^b$ and the observations $y^o$.  Then nudging method uses $K$ to modify the original equation in the following way:  
\begin{align}
\frac{\partial x^b(t,\bfx)}{\partial t} = M(t, x^b(t,\bfx)) + K(y^o(t) - h(x^b(t))),\label{eq: nudging in general}
\end{align}
in which $h(\cdot)$ is the observation operator. The flexibility of $K$ almost guarantees that there always exists such a $K$ that could adjust $x^b$ in the right direction. But the wealthiness of  degrees of freedom also brings difficulty to the determination of $K$. Thus explicit assumption on the structure of $K$ is necessary. For instance, the $K$ matrix  based on the Kalman gain matrix is used in \cite{Lei2011AHN}. In this manuscript, a displacement-based nudging strategy is considered in section \ref{subsection: displacement-based nudging}. The term $K(y^o-h(x^b))$ is chosen to be a discretised version of the Lie derivative $\mathcal{L}_{\bfu(t)}\theta_{x^b}(t)$, in which $\bfu(t)$ is a vector field calculated from $y^o(t)$ and $x^b(t)$, and $\theta_x^b(t)$ is a tensor field determined by the value of $x^b(t)$ and the governing equation of the state variable $x$. A brief introduction of the concept of Lie derivative and tensor fields is given in appendix A. The discretised $\mathcal{L}_{\bfu(t)}\theta_{x^b}$ can be written in the form $K(y^o(t) - h(x^b(t)))$, the explicit expression of which depends on how the model is discretised. Thus the displacement-based nudging strategy is indeed a particular case of Eq.\eqref{eq: nudging in general}. Since $K$ can be explicitly computed for each given $y^o(t)$ and $x^b(t)$, the displacement-based nudging strategy can also be applied to instantaneous data assimilation tasks like 3DVar and EnKF at a single time point. The instantaneous displacement-based nudging incorporated with the plain ensemble Kalman filter method, which we call morphed EnKF following \cite{Beezley2007MorphingEK}, is the method numerically tested in this manuscript.

To better understand our motivation, consider $x^b$ and $y^o$ as shown in Fig. \ref{fig: intro: xb and yo }. The difference simply relates to shifted positions of two bright spots. A natural way to nudge $x^b$ is the following:
\begin{figure}
	\centering
	\includegraphics[width = 0.9\textwidth]{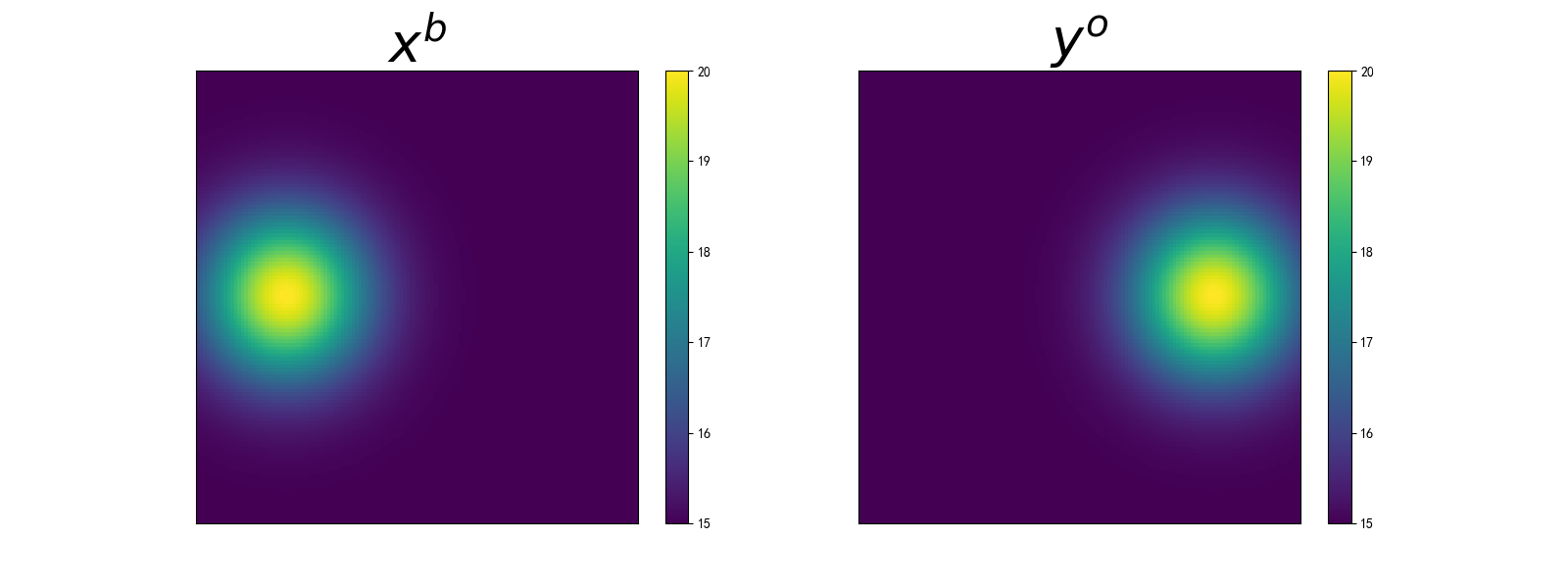}
	\caption{Suppose that the difference between the model estimate $x^b$ and observed state $y^o$ only differ by the position, then the most natural  nudging strategy  is to gradually moves $x^b$ rightward.}
	\label{fig: intro: xb and yo }
\end{figure}

\begin{itemize}
	\item[] (Displacement-based nudging) \hspace{3mm} $x^b(\bfx)\longleftarrow x^b((1-\epsilon){\bfx} + \epsilon T(\bfx))$. 
\end{itemize}
The displacement of the physical field is thus explained by the displacement of the points of the domain, which we denote by $T(\bfx)$. In Fig. \ref{fig: intro: xb and yo } the displacement map $T(\cdot)$ would map every point in the domain to a point on its left. It is determined by $y^o$ and $x^b$. $\epsilon$ is a parameter for nudging. As $\epsilon$ increases from $0$ to $1$, the displacement-based nudging strategy shall gradually shift the bright spot in $x^b$ to the correct position.  
For continuous time nudging, we can choose $\epsilon = \Delta t$, where $\Delta t$ is the model integration time step. For instantaneous data assimilation task, we can choose $\epsilon = 1$. It can be seen that the main feature of the physical field $x^b$ is preserved for both small and large $\epsilon$. Physically consistent estimation of $T(\cdot)$ for given background estimate $x^b$ and observation $y^o$ is the main objective of this manuscript.   Note that it is implicitly assumed in the example of Fig.\ref{fig: intro: xb and yo } that for a given displacement map $T(\cdot)$ and a given field $x^b$, the adjusted field should be the composition of $x^b$ and the $T$, denoted by $x^b\circ T$. However, this assumption is not necessary.

In general,  the position error of physical field can be represented by a diffeomorphism (an invertible smooth map whose inverse is also smooth)  $T:\Omega \longrightarrow\Omega$, where $\Omega$ refers to the physical domain under consideration. In the example of Fig.\ref{fig: intro: xb and yo }, $\Omega$ is the infinite $xy$ plane. We use $\mathcal{D}(\Omega)$ to denote the set of all diffeomorphisms of $\Omega$. Denote by $S_1$ the target field and $S_2$ the background field. In the scenario of data assimilation, $S_1$ refers to the observation field $y^o(t)$ and $S_2$ is the prior estimate of the same state variable $h(x^b(t))$.
A numerical evaluation of $T$ for two physical fields $S_1$ and $S_2$, can then target the minimization of a cost function of the following form:
\begin{align}
\begin{cases}
T = \displaystyle\argmin_{T\in \mathcal{D}(\Omega) } \|S_1 - T^\# S_2\|^2 + a\|T\|^2 \\
\text{Proper boundary conditions.}
\end{cases}\label{eq: global optimisation}
\end{align}
The $\#$ symbol here is not the musical sharp operator in differential geometry. $T^\#S_2$ represents the adjusted $S_2$ for a given displacement map $T$. The first term of \eqref{eq: global optimisation} thus evaluates the residue error of the adjusted $S_2$ for a specific $T$. It will be noted later that $T^{\#}S$ has been implicitly assumed, however questionably, to be $S\circ T$ or $S\circ T^{-1}$ by most of the literature in geophysics of the past. $\|T\|$, a chosen norm for the map $T$, serves as a regularization term. In some situation $\|T\|$ merely defines the complexity of $T$. But $\|T\|$ could be more complicated if the user has other needs. $a>0$ is a pre-chosen positive constant that weights this regularization term. Thus  the optimisation problem \eqref{eq: global optimisation} attempts to find a map $T(\cdot)$ that could reduce the residue error while $T(\cdot)$ itself is not too complex so that overfitting does not occur. Directly solving the optimisation problem \eqref{eq: global optimisation} may technically be  difficult. Instead, it might be easier to iteratively solve an "infinitesimal" version : 
\begin{align}
\begin{cases}
{\bf v}_i = \displaystyle\argmin_{{\bf v}\in T\Omega }\frac{\partial \|S_1 - {\bf \Phi}_{\bf v}(s)^\#S_2\|^2}{\partial s}\big{|}_{s=0} + a\|{\bf v}\|^2 \\
\text{ Proper boundary condition for ${\bf v}_i$},
\end{cases}
\label{eq: infinitesimal optimisation}
\end{align}
in which $i$ refers to the $i$-th iterative step, $T\Omega$ the set of tangent vector fields on $\Omega$. The ${\bf v}$'s are tangent vector fields on $\Omega$ with sufficient regularity and following some given boundary condition. The boundary condition of ${\bf v}$ usually requires ${\bf v}_i$ to have zero flux at $\partial \Omega$, the boundary of $\Omega$, i.e ${\bf v}_i\big{|}_{\partial \Omega}$ is always parallel to $\partial\Omega$. ${\bf \Phi}_{\bf v}(\cdot)$ is a flow of differentiable maps generated by ${\bf v}$:
\begin{align}
\begin{cases}
&\frac{\partial {\bf \Phi}_{\bf v}(s,\bfx)}{\partial s} = {
	\bf v}({\bf \Phi}_{\bf v}(s, \bfx))\text{, \hspace{2mm} for any $\bfx\in\Omega$ and $s\in [0, s_{\text{max}})$,}\\
&{\bf \Phi}_{\bf v}(0,\bfx) = \bfx, \hspace{2mm} \text{for any $\bfx\in\Omega$.}
\end{cases}\label{eq: T and flow}   
\end{align}
$s$ is a virtual time parameter in \eqref{eq: T and flow}. For each fixed $s$, ${\bf \Phi}_{\bf v}(s,\cdot): \Omega\to\Omega$ is a map from $\Omega$ to itself, which we write for simplicity as ${\bf \Phi}_{\bf v}(s)$.  For each fixed $\bfx_0\in\Omega$, the image of ${\bf \Phi}_{\bf v}(\cdot, \bfx_0)$, denoted by $\bfx(\cdot)$ in this paragraph, is a curve on $\Omega$ which satisfies the ordinary differential equation: $\frac{{d \bfx(s)}}{ds} = {\bf v} (\bfx(s))$ with initial condition $\bfx(0) = \bfx_0$. ${\bf \Phi}_{\bf v}$ would be an automorphism uniquely determined by \eqref{eq: T and flow} if ${\bf v}$ is Lipschitz continuous. ${\bf \Phi}_{\bf v}$ may not be unique (nor an automorphism) for merely continuous ${\bf v}$. However, for each point $\bfx\in\Omega$, different ${\bf \Phi}_{\bf v}$ would result in the same value of the first term in \eqref{eq: infinitesimal optimisation}.

The optimisation problem \eqref{eq: infinitesimal optimisation} is not derived from \eqref{eq: global optimisation}. While the first term in \eqref{eq: global optimisation} directly evaluates the residue error of the adjusted $S_2$, the first term in \eqref{eq: infinitesimal optimisation} measures how fast the residue error can be reduced if $S_2$ starts to morph along the vector field ${\bf v}$. The second term in \eqref{eq: infinitesimal optimisation} measures the complexity of ${\bf v}$. While \eqref{eq: global optimisation} directly defines the "optimal" $T$ using one equation, \eqref{eq: infinitesimal optimisation} aims to construct $T$ iteratively. At the $i$-th iterative step (denoted by the subscript $i$ in \eqref{eq: infinitesimal optimisation}), ${\bf v}_i$ is obtained by solving \eqref{eq: infinitesimal optimisation}. Then $S_2$ is adjusted along the direction of ${\bf v}_i$ for a small virtual time step $\epsilon$: $S_2\leftarrow T_i^{\#}S_2$ where $T_i = {\bf \Phi}_{{\bf v}_i}(\epsilon)$. The updated $S_2$ is then used as the input for the $(i+1)$-th iterative step to calculate ${\bf v}_{i+1}$. A stopping criterion is required to stop this iterative process. Then the final map $T$ can be constructed as: $T \approx T_N\circ T_{N-1}\circ \cdots T_0$. In many applications there is no need to calculate $T_i$ (hence $T$) explicitly. Only the final $T^{\#}S_2$ (hence all the $T_i^{\#}S_2$ in the middle steps) is needed and it could be derived efficiently based on the first order Taylor expansion of $T_i^{\#}S_2$. The specific definition of $T^{\#}S$ shall be given in section 2, together with the formula and some specific examples on how to calculate the first order Taylor expansion of ${\bf \Phi}_{{\bf v}_i}(\epsilon)^{\#}S$.

\begin{figure}
	\centering
	\includegraphics[width=\textwidth]{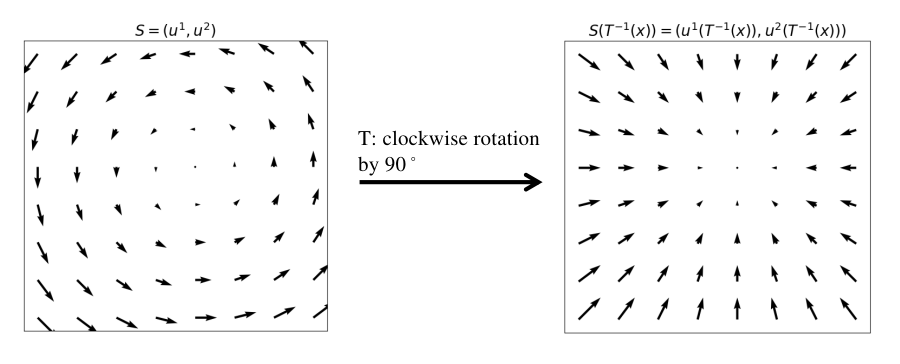}
	\caption{Suppose that the original field $S$ is a rotating wind field and the given displacement map $T$ is clockwise rotation by $90^\circ$. Then the direct composition of $S$ and $T$ results in a wind field of completely different feature.}
	\label{fig: 2}
\end{figure}

A specific optimisation problem in the form of \eqref{eq: global optimisation} (or \eqref{eq: infinitesimal optimisation}) is determined by two factors: the definition of $T^\#$ and the choice of the norms for the first and second terms in \eqref{eq: global optimisation} (or \eqref{eq: infinitesimal optimisation}). Intuitively, the definition of $T^{\#}$ is the rule of how to morph a physical field following a given map $T$.
Eq.\eqref{eq: infinitesimal optimisation}. The definition of $T^\#$ is essential and some relevant algorithms are reviewed or briefly discussed in section 2 and section 3, including the optical flow (OF) algorithms  \cite{HornSchunck1981, Younesetal2000, Lefvre2008OpticalFA}, the large deformation diffeomorphisms \cite{Trouv1995, Trouv1998DiffeomorphismsGA, Beg2005ComputingLD, Bruveris2009TheMM}, the metamorphoses  method \cite{Trouv2005LocalGO, Trouv2005MetamorphosesTL}, and some other algorithms \cite{Ravela2007DataAB, Beezley2007MorphingEK}. Although some of these algorithms have been applied to geophysical observations,  $T^{\#}S = S\circ T$ or $S\circ T^{-1}$ is commonly implicitly assumed, without any consideration on the dynamics of the geophysical fields. The displacement map $T$ may then cause severe structural errors  if it is used to adjust the displacement error of unobserved state variables. Fig.\ref{fig: 2} illustrates an example of directly applying $T$ by composition. The original physical field $S$ is a rotational vector field, rotating counter-clockwise. Suppose that a displacement map $T$ is given or   estimated from observations of other state variables. Suppose that the given $T$ is  a clockwise rotation by $90^\circ$. A direct application of $T$ then transforms $S$ to a displacement field $S\circ T^{-1}$ displaying completely different features.

To circumvent such undesirable results, \cite{Zhen2023PhysicallyCC} considered a differential geometry framework and the use of tensor fields to formulate $T^\#$ in Eq.\eqref{eq: global optimisation} and \eqref{eq: infinitesimal optimisation}.
In this new perspective, the choice of tensor fields can then be explicitly dictated by the dynamical equations of the underlying physical fields. $T^\#$ will follow dynamical principles, and certain physical quantities are naturally conserved during the morphing process. A similar argument has been given for small displacement cases \cite{Wildes2000RecoveringEO}. This further leads to a new alignment strategy in which the displacement flow calculated from the observed state variables can also be applied to partially correct the displacement of the physical fields of unobserved state variables. Such a new alignment strategy can naturally be incorporated with the classical ensemble Kalman filter to reduce inherent difficulties arising from linear algorithms.

%This manuscript is organized as follows. 
The physically consistent formulation of $T^\#$, based on the concept of tensor fields 
%in differential geometry, 
is given in section 2. This leads to a  version of the optimisation problem in the form of \eqref{eq: infinitesimal optimisation} different from the commonly used methods in geophysics \cite{HornSchunck1981, Ravela2007DataAB, Beezley2007MorphingEK}. Existence and uniqueness of the resulting solution are provided. A physical interpretation of the new optimisation problem is discussed,  with a brief review of some classical OF algorithms \cite{HornSchunck1981,Younesetal2000,Lefvre2008OpticalFA}. These OF algorithms are then compared with the proposed algorithm. The large deformation diffeomorphism strategy \cite{Beg2005ComputingLD,Bruveris2009TheMM,Trouv1995,Trouv1998DiffeomorphismsGA} and the metamorphoses strategy \cite{Trouv2005LocalGO,Trouv2005MetamorphosesTL} are also discussed in section 2. The new nudging strategy and its associated data assimilation strategy are presented in section 3. Differences between the proposed data assimilation strategy and some existing methods \cite{Ravela2007DataAB, Beezley2007MorphingEK} are  discussed. Using the thermal shallow water dynamical framework, numerical results are serving to support the proposed developments. Conclusion is given in Section 4. The complete code to reproduce the numerical results in this paper is available at http://doi.org/10.5281/zenodo.10252175.

\section{A differential geometry formulation of the optimisation problem}
A simple introduction of the concepts of differential geometry used in this manuscript, but restricted to the Euclidean case, is provided in appendix \ref{Appendix: differential geometry concepts}. For a thorough and mathematical introduction of tangent bundle, cotangent bundle, tensor fields, Riemannian manifold, the definition of ${\bf u}^{\musFlat{}}$, and other fundamental concepts in differential geometry, see for instance \cite{Chern1999LecturesDiffGeo}.
Let $(\Omega, g)$ be a compact oriented Riemannian manifold of dimension $n$ with  or without boundary, in which $g$ is the Riemannian metric. For any smooth vector field ${\bf u}$, ${\bf u}^{\musFlat{}}$ is the differential $1$-form such that ${\bf u}^{\musFlat{}}({\bf v}) = g({\bf u},{\bf v}) = \inp{{\bf u}}{{\bf v}}$ for any smooth vector field ${\bf v}$. Let $T\Omega$ be the tangent bundle of $\Omega$ and $T^*\Omega$ the cotangent bundle. For any diffeomorphism $\phi: \Omega\longrightarrow\Omega$, $\phi_*$ refers to the push-forward map of $T\Omega$ induced by $\phi$, and $\phi^*$ to the pull-back map of $T^*\Omega$ induced by $\phi$. For the simplicity of notation, we generalize the definition of push-forward map $\phi_*$ for diffeomorphism $\phi$, so that it also acts on $T^*\Omega$: $\phi_*\omega = (\phi^{-1})^*\omega$. Thus $\phi_*$ can act on both vector fields and differential forms at the same time. And for a differential form $\omega$, both $\phi^*$ and $\phi_*$ act on $\omega$, but in different ways. $\phi_*\omega = (\phi^*)^{-1}\omega$. An explanation of why this generalization is made is provided in appendix \ref{Appendix: explain T_* acts on differential forms}. It must be noted that the action of $\phi_*$ on $T^*\Omega$ is only restricted to this manuscript. This symbol is not used in this way in other literature.

$\phi_*$ further generalizes to an isomorphism of tensor fields of any specific type:
\begin{align}
\phi_*: V_1\otimes V_2\otimes\cdots\otimes V_l&\longrightarrow V_1\otimes V_2\otimes\cdots\otimes V_l \\
\alpha_1\otimes\alpha_2\otimes\cdots\otimes\alpha_l &\longrightarrow (\phi_*\alpha_1)\otimes\cdots\otimes(\phi_*\alpha_l)
\end{align}
in which $V_i$ can either be a copy of $T\Omega$ or a copy of $T^*\Omega$.

For the smooth vector field ${\bf u}$ such that ${\bf u}\big{|}_{\partial\Omega}\in i_*(T\partial\Omega)$, where $i:\partial\Omega\longrightarrow\Omega$ is the natural embedding map of the boundary and $i_*$ is the push-forward map, induced by $i$, that acts on the vector fields on $\partial \Omega$.   ${\bf \Phi}_{\bf u}: [0,\epsilon]\times\Omega\longrightarrow\Omega$ denotes the flow that satisfies
\begin{align}
\frac{\partial{\bf \Phi}_{\bf u}(s,{\bf x})}{\partial s} = {\bf u}({\bf \Phi}_{\bf u}(s,{\bf x})).
\end{align}
Using our notation ${\bf \Phi}_{\bf u}(s)_*$, the Lie derivative of  a tensor field $\theta$ with respect to $u$ is usually defined as:
\begin{align}
\mathcal{L}_{\bf u}\theta = \lim_{s\to 0}\frac{{\bf \Phi}_{-{\bf u}}(s)_*\theta - \theta}{s},
\end{align}
which is also a tensor field, of the same type as $\theta$. To first order, a Taylor expansion gives: 
\begin{align}
{\bf \Phi}_{\bf u}(s)_*\theta = \theta - s\mathcal{L}_{\bf u}\theta + o(s).\label{eq: Tayler expansion for Phiutheta}
\end{align}

A brief introduction of the operators mentioned in this paragraph can be found in appendix \ref{Appendix: differential geometry concepts}. The Riemannian metric $g$ provides the inner product for the tangent space $T_{\bf x}\Omega$ for any ${\bf x}\in\Omega$. This inner product could be generalized to the space of tensor fields of any specific type. For $\theta_1,\theta_2\in V_1\otimes\cdots\otimes V_k$, the inner product induced by $g$ at ${\bf x}\in\Omega$ is written as $\inp{\theta_1}{\theta_2}_{\bf x}$. The corresponding norm is written as $|\theta_2|_{\bf x}^2 = \inp{\theta_2}{\theta_2}_{\bf x}$. The Riemannian metric further induces several operators. For any integer $0\leq k\leq n$, let $\Lambda^k\Omega$ be the space of smooth differential $k-$forms. Let $*: \Lambda^{k}\Omega\to\Lambda^{n-k}\Omega$ be the Hodge star operator, $d:\Lambda^k\Omega\to\Lambda^{k+1}\Omega$ the exterior derivative, and $\delta = (-1)^{n(k-1)+1}*d*:\Lambda^{k+1}\Omega\to\Lambda^k\Omega$ the codifferential operator.  The Hodge-Laplacian operator is defined as $\Delta = d\delta + \delta d: \Lambda^k\Omega\to\Lambda^k\Omega$.  Let $dV$ be the volume form on $\Omega$.

Recall that ${\bf u}^{\musFlat{}}$ is the differential $1$-form determined by vector field ${\bf u}$ and the Riemannian metric tensor $g$. With the above notations, we now provide a specific formulation of the optimisation problem \eqref{eq: infinitesimal optimisation} for tensor fields on a Riemannian manifold:
\begin{definition}
	Let $\theta_1$ and $\theta_2$ be two tensor fields on $\Omega$ of the same type. When $\Omega$ is a compact oriented Riemannian manifold without boundary, solve for
	\begin{align}
	{\bf u} = \argmin_{{\bf u}\in T\Omega} \int_{\Omega}\bigg{[}\frac{\partial |\theta_1 - {\bf \Phi}_{\bf u}(s)_*\theta_2|_{\bf x}^2}{\partial s}\Big{|}_{s=0} + a_1|d{\bf u}^{\musFlat}|_{\bf x}^2 + a_1|\delta {\bf u}^{\musFlat{}}|_{\bf x}^2 + a_0|{\bf u}^{\musFlat{}}|_{\bf x}^2\bigg{]} dV({\bf x}).\label{eq: optical flow for tensor fields, compact no boundary}
	\end{align}
	When $\Omega$ is a compact oriented Riemannian manifold with boundary, solve for
	\begin{align}
	\begin{cases}
	{\bf u} = \displaystyle\argmin_{{\bf u}\in T\Omega} \int_{\Omega}\bigg{[}\frac{\partial |\theta_1 - {\bf \Phi}_{\bf u}(s)_*\theta_2|_{\bf x}^2}{\partial s}\Big{|}_{s=0} + a_1|d{\bf u}^{\musFlat{}}|_{\bf x}^2 + a_1|\delta {\bf u}^{\musFlat{}}|_{\bf x}^2 \bigg{]} dV({\bf x})\\
	{\bf u}\big{|}_{\partial\Omega} = i_*{\bf v},
	\end{cases}\label{eq: optical flow for tensor fields, compact with boundary}
	\end{align}
	where ${\bf v}\in T\partial\Omega$ is the given boundary condition and $i: T\partial\Omega\to\Omega$ is the natural embedding of $\partial\Omega$.
\end{definition}

Using the Taylor expansion of ${\bf \Phi}_{\bf u}(s)_* \theta_2$, see Eq.\eqref{eq: Tayler expansion for Phiutheta}, Eq.\eqref{eq: optical flow for tensor fields, compact no boundary} is equivalent to 
\begin{align}
{\bf u} = \argmin_{{\bf u}\in T\Omega} \int_{\Omega}\bigg{[}2\inp{\theta_1-\theta_2}{\mathcal{L}_{\bf u}\theta_2}_{\bf x} + a_1|d{\bf u}^{\musFlat{}}|_{\bf x}^2 + a_1|\delta {\bf u}^{\musFlat{}}|_{\bf x}^2 + a_0|{\bf u}^{\musFlat{}}|_{\bf x}^2\bigg{]} dV({\bf x}),\label{eq: optical flow for tensor fields, compact no boundary, v2}
\end{align}
and Eq.\eqref{eq: optical flow for tensor fields, compact with boundary} equivalent to 
\begin{align}
\begin{cases}
{\bf u} = \displaystyle\argmin_{{\bf u}\in T\Omega} \int_{\Omega}\bigg{[}2\inp{\theta_1-\theta_2}{\mathcal{L}_{\bf u}\theta_2}_{\bf x} + a_1|d{\bf u}^{\musFlat{}}|_{\bf x}^2 + a_1|\delta {\bf u}^{\musFlat{}}|_{\bf x}^2 \bigg{]} dV({\bf x})\\
{\bf u}\big{|}_{\partial\Omega} = i_*{\bf v}.
\end{cases}\label{eq: optical flow for tensor fields, compact with boundary, v2}
\end{align}

\begin{theorem}\label{thm: existence and uniqueness for tensor fields}
	For smooth tensor fields and vector fields $\theta_1,\theta_2,$ and ${\bf v}$, the optimisation problems \eqref{eq: optical flow for tensor fields, compact no boundary, v2} and \eqref{eq: optical flow for tensor fields, compact with boundary, v2} are always uniquely solvable. And the solution is smooth.
\end{theorem}
When $\Omega$ is an Euclidean manifold, an explicit  guidance on how to calculate everything, including the exterior derivative operator, the codiffenrential operator, the Hodge Laplacian, and the Lie derivative, is provided in appendix \ref{Appendix: differential geometry concepts}. Theorem \ref{thm: existence and uniqueness for tensor fields} is a direct consequence of some proven mathematical results in \cite{Duff1952HARMONICTO} (or \cite{Schwarz1995HodgeD}). A complete demonstration is provided in the appendix. The formulation of \eqref{eq: optical flow for tensor fields, compact no boundary} and \eqref{eq: optical flow for tensor fields, compact with boundary}, together with theorem \eqref{thm: existence and uniqueness for tensor fields} provide an option for the theory of defining and computing the displacement flow of two tensor fields. Due to Poincare lemma, both the regularization terms in \eqref{eq: optical flow for tensor fields, compact no boundary} and \eqref{eq: optical flow for tensor fields, compact with boundary} are equivalent to the $H^1$ norm. Still, the regularization term does not have to be the $H^1$ norm. For instance, $|d\delta{\bf u}^{\musFlat{}}|_{\bf x}^2 + |\delta d {\bf u}^{\musFlat{}}|_{\bf x}^2$ is considered in \cite{Corpetti2002DenseEO}. However, the vector field ${\bf u}$ calculated from Eq.\eqref{eq: optical flow for tensor fields, compact no boundary} or \eqref{eq: optical flow for tensor fields, compact with boundary}
should not necessarily be the same as the true velocity field, largely depending on the regularization terms.
%This comes from the regularization terms in Eq.\eqref{eq: optical flow for tensor fields, compact no boundary} and \eqref{eq: optical flow for tensor fields, compact with boundary}, merely constructed for mathematical reasons. 
True physical laws may not always be well represented by such regularization terms.

\subsection{A detailed formulation of $T^\#$ 
	%and  \eqref{eq: infinitesimal optimisation} 
	for physical fields}
We can now write down a more explicit expression for \eqref{eq: infinitesimal optimisation}. Given two physical fields $S_1({\bf x})$ and $S_2({\bf x})$ that represent the same  state variable $S$, we first associate a tensor to this state variable $S$. Thus $S_1$ and $S_2$ corresponds to two tensor fields, $\theta_1$ and $\theta_2$, of the same type. The displacement vector field ${\bf u}$ is then calculated based on \eqref{eq: optical flow for tensor fields, compact no boundary} or \eqref{eq: optical flow for tensor fields, compact with boundary} for $\theta_1$ and $\theta_2$. It is required that the correspondence between $S$ and the tensor is 1-1 and onto, i.e.  a unique value of $S$ can be inferred from a given value of the tensor, and vice versa. Thus, for a given displacement map (a diffeomorphism of $\Omega$) $T:\Omega\longrightarrow\Omega$, $T^\#S_2$ is defined in three steps:
\begin{itemize}
	\item[] 1) construct tensor field $\theta_2$ based on $S_2$;
	\item[] 2) $\theta_2^{\text{new}}\leftarrow T_*\theta_2$;
	\item[] 3) $T^\#S_2$ is then inferred from $\theta_2^{\text{new}}$.
\end{itemize}

Under $T^\#$, certain physical quantities will naturally be conserved as long as $T$ is a diffeomorphism of $\Omega$ which preserves the orientation. Iteratively solving \eqref{eq: infinitesimal optimisation}, a final displacement  map $T$ can then be constructed from $T_i={\bf \Phi}_{{\bf u}_i}(\epsilon)$: $T = T_{N-1}\circ \cdots\circ T_{0}$, where $\epsilon$ is the pre-chosen time step for the virtual time variable $s$. We demonstrate the conservative nature of $T^\#$ using the following examples.
\begin{example}\label{example: density}
	Suppose $S$ denotes the density of a flow. We can associate $S$ to a differential $n-$form: $\theta_S = SdV$, where $dV$ is the volume form on $\Omega$. If $\Omega$ itself is $n$-dimensional and is part of the $n$-dimensional Euclidean space, then $dV = dx_1\wedge \cdots\wedge dx_n$. Then for any diffeomorphism $T$ that preserves the orientation of $\Omega$, 
	\begin{align}
	(T_*\theta_S)({\bf x}) = ((T^{-1})^*\theta_S)({\bf x}) = S(T^{-1}({\bf x})) (T^{-1})^*(dV) = S(T^{-1}({\bf x}))J_{T^{-1}}({\bf x}) dV(\bfx),
	\end{align} 
	in which $J_{T^{-1}}$ is the Jacobian determinant of $T^{-1}$. Recall that $dV$ is a differential $n$-form and $(T^{-1})^*dV$ is the pull-back of $dV$.  Thus for a given displacement map $T$ and a given density field $S$, under the correspondence that $S\longleftrightarrow SdV, $ the adjusted field at ${\bf x}\in \Omega$ has value  $(T^\#S)({\bf x}) = S(T^{-1}({\bf x}))J_{T^{-1}}({\bf x})$. If $\Omega$ itself is $n$-dimensional and is part of a $n$-dimensional Euclidean space and suppose that 
	\begin{align}
	T^{-1}(x_1,...,x_n) = ((T^{-1})_1(x_1,...,x_n), (T^{-1})_2(x_1,...,x_n),\cdots, (T^{-1})_n(x_1,...,x_n)) 
	\end{align}
	then 
	\begin{align}
	T^{\#}S = S(T^{-1}({\bf x}))\det \begin{pmatrix}
	\frac{\partial (T^{-1})_i}{\partial x_j}
	\end{pmatrix}_{1\leq i,j\leq n}
	\end{align}
	For pull-back map $(T^{-1})^*$, it yields that $\displaystyle\int_{\Omega}SdV = \int_{\Omega}(T^{-1})^*(SdV)$. Thus 
	\begin{align}
	\int_{\Omega}SdV = \int_{\Omega}T^\#S dV,
	\end{align}
	which means that the total mass is conserved.
\end{example}

\begin{example}\label{example: vorticity}
	Suppose that $S = {\bf u} = (u_1,u_2)$ is a velocity field on a 2 dimensional domain in $\mathbb{R}^2$. We can choose to associate to $S$ a differential $1-$form $\theta_S = u_1 dx_1 + u_2 dx_2$. Then for any 1-1 and onto map $T$ which preserves the orientation of $\Omega$, direct calculation of pull-back operator yields
	\begin{align}
	&((T^{-1})^*\theta_S)({\bf x})\nonumber \\
	=& \Big{[}u_1(T^{-1}({\bf x}))\frac{\partial (T^{-1})_1}{\partial x_1} + u_2(T^{-1}({\bf x}))\frac{\partial (T^{-1})_2}{\partial x_1}\Big{]} dx_1 + \Big{[}u_1(T^{-1}({\bf x}))\frac{\partial (T^{-1})_1}{\partial x_2} + u_2(T^{-1}({\bf x}))\frac{\partial (T^{-1})_2}{\partial x_2}\Big{]} dx_2
	\end{align}
	This shows that, under the correspondence $(u_1,u_2)\longleftrightarrow u_1dx_1+u_2dx_2$, for a given velocity field $(u_1,u_2)$ and a given displacement map $T$, the adjusted velocity field $T^\# (u_1,u_2)$ is 
	\begin{align}
	\bigg{(}u_1(T^{-1}(\bfx))\frac{\partial (T^{-1})_1}{\partial x_1} + u_2(T^{-1}(\bfx))\frac{\partial (T^{-1})_2}{\partial x_1},u_1(T^{-1}(\bfx))\frac{\partial (T^{-1})_1}{\partial x_2} + u_2(T^{-1}(\bfx))\frac{\partial (T^{-1})_2}{\partial x_2}\bigg{)}
	\end{align}
	Since the exterior derivative operator $d$ commutes with the pull-back operator $d(T^{-1})^* = (T^{-1})^*d$, 
	we have
	\begin{align}
	\int_{\Omega}d\theta_S = \int_{\Omega}(T^{-1})^*(d\theta_S) = \int_{\Omega}d((T^{-1})^*\theta_S) = \int_{\Omega}d(T^{\#}\theta_S)
	\end{align}
	Note that $d\theta_S = (\frac{\partial u_2}{\partial x_1} - \frac{\partial u_1}{\partial x_2})dx_1\wedge dx_2$. Therefore
	\begin{align}
	\int_{\Omega}(\frac{\partial u_2}{\partial x_1} - \frac{\partial u_1}{\partial x_2})dx_1\wedge dx_2 = \int_{\Omega}(\frac{\partial [T^\#(u_1,u_2)]_2}{\partial x_1} - \frac{\partial [T^\#(u_1,u_2)]_1}{\partial x_2})dx_1\wedge dx_2,
	\end{align}
	in which $[T^\#(u_1,u_2)]_i$ refers to the $i-$th component of $[T^\#(u_1,u_2)]$ for $i=1,2$.
	Hence the total vorticity is conserved.
\end{example}

\subsection{A physical interpretation of the tensor field definitions}
Described above for \ref{example: density} and \ref{example: vorticity}, the conservative nature of $T^\#$ largely depends on the choice of tensor fields. The question now translates to how to choose the adequate tensor fields. Since Eq.\eqref{eq: Tayler expansion for Phiutheta} describes how the tensor field changes along the vector field ${\bf u}$, the Lie derivative in \eqref{eq: Tayler expansion for Phiutheta} shall closely relate to the transport equation of physical fields. In examples (\ref{example: transport eqn, 0-form}) and (\ref{example: transport eqn, n-form}), a flow $\bfPhi$ is constructed based on the dynamical equation of the underlying state variable $S$. It is then explained that the tensor field $\theta_S$ should be chosen in a way so that, when there is no external forcing or sources, $\theta_S$ instead of S, is conserved along the flow trajectories: i.e $\frac{d}{dt}\bfPhi(t)_*\theta_S(t,\bfx) = 0$. As a consequence, the Lie derivative $\mathcal{L}_\bfv\theta_S$ should exactly match the transportation term in the original dynamical equations. The correspondence between $S$ and $\theta_S$ then induces the explicit definition of $T^\#S$ for any diffeomorphism $T$. Note that the flow map $\bfPhi$ constructed in this subsection is completely determined by the original dynamical equation of $S$. It does not refer to the aforementioned displacement flow constructed for two physical fields $S_1$ and $S_2$. When $S$ refers to the velocity field, the corresponding tensor field $\theta_S$ can be chosen based on the principle of what velocity-related quantity should be conserved during the morphing process. 

%See the following two examples.
\begin{example}\label{example: transport eqn, 0-form}
	Consider $\Omega \subset \mathbb{R}^2$ is a two dimensional domain, and the original dynamical equation for $S$  
	\begin{align}
	\frac{dS}{dt} := S_t + {\bf v}\cdot\nabla S = \frac{\partial S}{\partial t} + v_1\frac{\partial S}{\partial x_1} + v_2\frac{\partial S}{\partial x_2} = 0,\label{eq: transport of point value}
	\end{align}
	in which ${\bf v}(t,{\bf x}) = (v_1(t,{\bf x}), v_2(t,{\bf x}))$ is the true velocity field in the system. Define the flow ${\bf \Phi}(t,{\bf x})$:
	\begin{align}
	\begin{cases}
	\frac{\partial {\bf \Phi}}{\partial t}(t,\bfx) = {\bf v}(t,{\bf \Phi}(t,\bfx)) \\
	{\bf \Phi}(0,\bfx) = \bfx.
	\end{cases}
	\end{align}
	
	Or in index notation, the above equation can be written as
	\begin{align}
	\begin{cases}
	\frac{\partial \Phi_1}{\partial t}(t,x_1,x_2) = v_1(t,\Phi_1(t,x_1,x_2),\Phi_2(t,x_1,x_2))\\
	\frac{\partial \Phi_2}{\partial t}(t,x_1,x_2) = v_2(t, \Phi_1(t,x_1,x_2), \Phi_2(t,x_1,x_2))\\
	\Phi_1(0,x_1,x_2) = x_1\\
	\Phi_2(0,x_1,x_2) = x_2,
	\end{cases}
	\end{align}
	where $\bfx=(x_1,x_2)$, ${\bf \Phi}(t,\bfx) = (\Phi_1(t,x_1,x_2), \Phi_2(t,x_1,x_2))$.
	
	Then Eq.\eqref{eq: transport of point value} means that $\frac{d}{dt}S(t, {\bf \Phi}(t,\bfx)) = 0$. Thus $S(t, {\bf \Phi}(t,\bfx)) = S(0,\bfx)$ for any $t > 0$, i.e.  $S$ is conserved along the trajectory of each point. In this case, ${\bf \Phi}(t): \Omega\longrightarrow \Omega$ plays the role of a ``true" displacement map. Recall that the goal is to find a displacement map $T$ for given snapshots $S_1$ and $S_2$. Choose the tensor $\theta = S$, which is simply a function (or differential 0-form as a tensor field), the relation $S(t, {\bf \Phi}(t,\bfx)) = S(0,\bfx)$  can be written as $\theta(t) = {\bf \Phi}(t)_*\theta(0)$. Correspondingly, we have $\mathcal{L}_{\bf u}\theta = {\bf u}\cdot\nabla S = u_1\frac{\partial S}{\partial x_1} + u_2\frac{\partial S}{\partial x_2}$.
	
	If the original dynamical equation for the physical field $S$ is modified 
	\begin{align}
	\frac{dS}{dt} := S_t + \bfv\cdot \nabla S = \text{(forcing terms)}
	\end{align}
	$S(0,\bfx)\neq S(t, {\bf \Phi}(t,\bfx))$. However, we can still choose $\theta = S$. The minimization problem \eqref{eq: infinitesimal optimisation} shall search for the vector field ${\bf u}$ transporting $S_2$ towards $S_1$ the fastest along the virtual flow ${\bf \Phi}_{\bf u}$. It is also equivalent to look for a virtual vector field ${\bf u}$ so that if $S_2$ is transported along ${\bf u}$ by the virtual system
	\begin{align}
	\begin{cases}
	\frac{\partial\tilde{S}}{\partial s} + {\bf u}\cdot\nabla\tilde{S} = 0\\
	\tilde{S}(0,\bfx) = S_2
	\end{cases}
	\end{align}
	in a virtual time interval $[0,\epsilon]$, then $S_2$ transforms the most rapidly towards $S_1$. 
\end{example}

\begin{example}\label{example: transport eqn, n-form}
	Now, consider $\Omega\subset\mathbb{R}^n$ an $n$-dimensional domain with the original dynamical equation  
	\begin{align}
	S_t + \nabla\cdot(S\bfv) = 0,  \label{eq: transport of density}
	\end{align}
	with $\bfv(t,\bfx)$ the true velocity field. Similar to the previous example, define the flow $\bfPhi(t,\bfx)$:
	\begin{align}
	\begin{cases}
	\frac{\partial\bfPhi}{\partial t } (t,\bfx) = \bfv(t,\Phi(t,\bfx))\\
	\bfPhi(0,\bfx) = \bfx
	\end{cases}
	\end{align}
	and the differential $n$-form $\theta(t,\bfx) = S(t,\bfx) dx_1\wedge dx_2\wedge\cdots\wedge dx_n$. Distinct from example \ref{example: transport eqn, 0-form}, Eq.\eqref{eq: transport of density} implies that $\bfPhi(t)^*\theta(t,\bfx) = \theta(0,\bfx)$, or equivalently, $\theta(t,\bfx) = (\bfPhi(t)^{-1})^*\theta(0,\bfx) = \bfPhi(t)_*\theta(0,\bfx)$. Indeed, direct calculation yields:
	\begin{align}
	&\bfPhi(t+\Delta t)^*\theta(t+\Delta t,\bfx) - \bfPhi(t)^*\theta(t,\bfx)\\
	=& 
	% (\Phi_{t\to t+\Delta t}\circ\Phi(t))^*\theta(t+\Delta t)
	\bfPhi(t)^* \bfPhi_{t\to t+\Delta t}^* \theta(t+\Delta t)  
	% - \Phi(t)^*\theta(t+\Delta t) + \Phi(t)^*(\theta(t+\Delta t) - \theta(t))\\
	- \bfPhi(t)^* \theta(t)\\
	=& \bfPhi(t)^*\Big{[}
	% \Phi_{t\to t+\Delta t}^*\theta(t+\Delta t) - \theta(t+\Delta t) 
	(\bfPhi_{t\to t+\Delta t}^* - Id) \theta(t+\Delta t)  
	+ \theta(t+\Delta t) - \theta(t)\Big{]}\\
	=& \bfPhi(t)^*\Big{[}
	% \Phi_{t\to t+\Delta t}^*\theta(t) - \theta(t) 
	(\bfPhi_{t\to t+\Delta t}^* - Id) \theta(t) 
	+ \theta(t+\Delta t) - \theta(t)  \Big{]} + o(\Delta t)\\
	=& \Delta t\bfPhi(t)^*\Big{[} \mathcal{L}_{\bfv(t)}\theta(t) + \theta_t(t)\Big{]} + o(\Delta t)\\
	=&\Delta t\bfPhi(t)^*\Big{[} \Big{(}S_t(t,\bfx) + \nabla\cdot(S(t,\bfx)\bfv(t))\Big{)}dx^1\wedge\cdots\wedge dx^n\Big{]} + o(\Delta t)\\
	=& o(\Delta t)
	\end{align}
	It shows that $\bfPhi(t)^*\theta(t,\bfx)$ does not change with time. 
\end{example}
In both examples, the choice of tensor $\theta$ is determined by the original transport equation \eqref{eq: transport of point value} or \eqref{eq: transport of density}. This choice then implies $T^\#$, and the minimization problem \eqref{eq: infinitesimal optimisation} will search for a vector field $\bfu$ that most rapidly transports $S_2$ towards $S_1$ by a virtual flow $\bfPhi_\bfu$. For geophysical fields, $T^\#$ is thus essential for the application of this method, and further discussed in section 3.

\subsection{Comparison with OF algorithms}
Given a time series of snapshots (i.e. brightness fields in $\Omega\subset\mathbb{R}^2$) $\{...,S(t-\Delta t), S(t), S(t+\Delta t),...\}$, the \cite{HornSchunck1981} OF algorithm aims at recovering the true velocity field $\bfu(t,\bfx)$  by minimizing the following functional (with proper boundary conditions):
\begin{align}
\bfu(t) = \arg\min_{\bfu} \|S_t + \inp{\nabla S}{\bfu}\|^2 + a_1\|\nabla \bfu\|^2 = \argmin_{\bfu}\int_{\Omega} [S_t + \inp{\nabla S}{\bfu}]^2 + a_1|\nabla \bfu|^2 d^n\bfx. \label{eq: optical flow}
\end{align}

We remark that the first term implies that the hidden $T^\#$ in \eqref{eq: optical flow} is simply $T^\#S = S\circ T^{-1}$ for the brightness field $S$. Or equivalently, $T^\# S$ is inferred from $T_*\theta_S$, where $\theta_S = S$ is a differential $0-$form. Indeed, for any vector field $\bfu$,
let $\bfPhi_\bfu(s,\bfx): [0,\Delta t]\times \Omega \longrightarrow \Omega$ be the following flow of the points in the domain:
\begin{align}
\begin{cases}
\frac{\partial \bfPhi_\bfu(s,\bfx)}{\partial s} = \bfu(\bfPhi_\bfu(s,\bfx))\\
\bfPhi_\bfu(0,\bfx) = \bfx.
\end{cases}
\end{align}

The first term in \eqref{eq: optical flow} represents the material derivative of $S(t+s,\bfPhi_\bfu(s, \bfx))$ with respect to $s$ at $s = 0$. The optimisation problem \eqref{eq: optical flow} searches for vector field ${\bf u}$ so that $S(t,\bfx) \approx S(t+\Delta t, \bfPhi_\bfu(\Delta t,\bfx))$. Let $S_2(\bfx) = S(t,\bfx)$, $S_1(\bfx) = S(t+\Delta t,\bfx)$, and $T(\bfx) = \bfPhi_\bfu(\Delta t,\bfx)$. Hence the first term in \eqref{eq: optical flow} searches for the vector field $\bfu$ so that $S_1(\bfx) \approx S_2(T^{-1}(\bfx))$. Clearly, this implies that $T^\# S_2 = S_2\circ T^{-1}$. \cite{Lefvre2008OpticalFA} generalizes this original OF method to the Riemannian manifold context. Still, the formulation of the problem in \cite{Lefvre2008OpticalFA}
implies that $T^\#S_2 = S_2\circ T^{-1}$.

In \cite{Younesetal2000}, an other variant is proposed, to minimize
\begin{align}
\bfu(t) = \argmin_\bfu \|S_t + \nabla\cdot(S\bfu)\|^2 + a_1\|\nabla \bfu\|^2 = \argmin_{\bfu}\int_{\Omega} [S_t + \nabla \cdot(S\bfu)]^2 + a_1|\nabla \bfu|^2 d^n\bfx. \label{eq: optical flow variant 1: conserve total brightness}
\end{align}
Let $S_1(\bfx) = S(t+\Delta t,\bfx)$ and $S_2(\bfx) = S(t,\bfx)$. Following the analysis in example \ref{example: transport eqn, n-form}, Eq.\eqref{eq: optical flow variant 1: conserve total brightness} aims at finding the vector field $\bfu$ so that $\bfPhi_\bfu^\#S_2\approx S_1$, in which $\bfPhi_\bfu^\#S$ is inferred from $\bfPhi_\bfu^\#\theta_S$ for $\theta_S = Sdx^1\wedge\cdots\wedge dx^n$. Thus, it can be concluded that the method of \cite{Younesetal2000} implies that $T^\#$ should be defined as if the physical fields are associated to differential $n$-forms.  
It has been explicitly stated in \cite{Wildes2000RecoveringEO} that the first residue error term in the formulation of optical flow method should depend on the physics law of the underlying state variable.
%. he choice of data term in OF-like algorithms should depend on the physic's law. This 
This corresponds to physically consistent choices of tensor fields in the present framework. The proper choices of tensor fields shall indeed also enable us to transport unobserved state variables, consistent with the underlying dynamics. OF algorithms can thus be reformulated for snapshots corresponding to snapshots of tensor fields.
%It is then not hard to generalize the optical flow algorithms to the case when the snapshots are snapshots of tensor fields.
For a time series $\{\cdots, \theta(t-\Delta t), \theta(t), \theta(t+\Delta t),\cdots\}$ of snapshots of tensor fields, the generalized OF method is to find the vector field $\bfu$:
\begin{align}
\bfu(t) = \argmin_{\bfu}\int_{\Omega}|\theta_t - \mathcal{L}_\bfu\theta|_\bfx^2 + a_1|\nabla \bfu|_\bfx^2 dV(\bfx),
\end{align}
in which $\bfu$ satisfies proper boundary conditions. Here we provide a theorem for the existence of the OF solution for Dirichlet boundary condition. Note that in  \cite{HornSchunck1981}, $\bfu$ is only required to have zero normal component at the boundary. 

\begin{theorem}\label{thm: OF of tensor fields}
	For fixed $t$, assume that both tensor fields $\theta$ and $\theta_t$ and the vector field $\bfv$ have finite $H^1$ norm. Further assume that $|d\theta|_\bfx,|\delta\theta|_\bfx$ are bounded in $\Omega$. Then for compact oriented Riemannian manifold $\Omega$ without boundary, the following minimisation problem has a unique solution which has a finite $H^1$ norm:
	\begin{align}
	\bfu(t) = \argmin_{\bfu\in H^1(T\Omega)} \int_{\Omega}|\theta_t - \mathcal{L}_\bfu\theta|_\bfx^2 + a_0|\uflat|_\bfx^2 + a_1|d\uflat|_\bfx^2 + a_1|\delta\uflat|_\bfx^2 dV(\bfx).\label{eq: OF of tensor fields, no boundary}
	\end{align}
	For a compact oriented Riemannian manifold with boundary, the following minimisation problem has a unique solution which has a finite $H^1$ norm:
	\begin{align}
	\begin{cases}
	\bfu(t) = \displaystyle\argmin_{\bfu\in H^1(T\Omega)} \int_{\Omega}|\theta_t - \mathcal{L}_\bfu\theta|_\bfx^2 + a_1|d\uflat|_\bfx^2 + a_1|\delta\uflat|_\bfx^2 dV(\bfx)\\
	\bfu(t)\big{|}_{\partial\Omega} = \bfv.
	\end{cases}\label{eq: OF of tensor fields, with boundary}
	\end{align}
\end{theorem}
The existence, uniqueness, and smoothness of order 1 of the solution is a direct consequence of Riesz representation (or Lax-Milgram) theorem. A mathematical proof is provided in the appendix.

While the proposed algorithm can be applied whenever two snapshots are provided, the OF algorithms require that the time between two consecutive images to be small, to best approximate $\theta_t$. To overcome such a constrain, the standard Horn and Schunck algorithm can be incorporated within a multi-scale strategy \cite{Meinhardt2013HornSchunckOF}, and wavelets \cite{Drian2013WaveletsAO}. Note, an infinitesimal formulation of the Horn and Schunck algorithm in the Euclidean space is proposed in \cite{Drian2013WaveletsAO}. The (generalized) OF method is closely related to the metamorphoses method \cite{Trouv2005MetamorphosesTL}, which is discussed in the next subsection.

\subsection{Comparison with other methods}
\cite{Trouv2005MetamorphosesTL} proposes a variational method to find the optimal way to transport a point $S_2\in\mathcal{M}$ to another point $S_1\in\mathcal{M}$ at a minimal cost. Here $\mathcal{M}$ is a Riemannian manifold. For geophysical applications, this manifold $\mathcal{M}$ can be interpreted as the phase space. Hence each point $S_i\in\mathcal{M}$ represents a state vector. It can be interpreted in a way that, to transport $S_2$ to $S_1$, \cite{Trouv2005MetamorphosesTL} does not only use vector fields, but also an external forcing to influence the state at each time step. The solution of the variational problem then consists of a time sequence of vector fields and a time sequence of "external forcing". Considering an "external forcing", the eventual state $T^\#S_2$ can exactly match the target state $S_1$. This method can be applied to general $T^{\#}$ and follows from rigorous mathematical developments \cite{Trouv2005LocalGO}, but turns out to be very computationally demanding. Interestingly, the metamorphoses method is closely related to the optical flow algorithm. Suppose a high frequency time series of snapshots is given: $S_n = S(n\Delta t)$, for $n = 0,1,2,...$. Let $\theta_n$ be the tensor field corresponding to $S_n$. If metamorphoses method is applied to estimate a vector field $\bfu_n$ to transport $\theta_{n}$ to $\theta_{n+1}$, then the cost function of metamorphoses method can be approximately written as
\begin{align}
\bfu_n = \argmin_{\bfu\in T\Omega}\int_{\Omega}\Big{|} \frac{\theta_{n+1} - (\theta_n - \Delta t\mathcal{L}_\bfu\theta_n) }{\Delta t}\Big{|}_\bfx^2 + |\bfu|_\bfx^2 dV(\bfx).\label{eq: metamorphoses for high frequency data}
\end{align}
The quantity inside the first norm is the external forcing that metamorphoses method considers. A vector field $\bfu$ transports $\theta_n$ to $\theta_n^+ = \theta_n-\Delta t\mathcal{L}_{\bfu}\theta_n$. Then an additional forcing $F_n(\bfu) = \theta_{n+1} - \theta_n^+$ needs to be added to $\theta_n^+$ so that the final tensor field exactly is $\theta_{n+1}$. The first term in \eqref{eq: metamorphoses for high frequency data} thus evaluates the price of the external forcing $F_n(\bfu)$.  The second term is again the smoothness constraint of $\bfu$. It is  not hard to see that in this case the metamorphoses method \eqref{eq: metamorphoses for high frequency data} minimizes approximately the same cost function as the (generalized) optical flow method \eqref{eq: OF of tensor fields, no boundary}-\eqref{eq: OF of tensor fields, with boundary}.

Some efforts about large deformation diffeomorphic metric matching (LDDMM) \cite{Trouv1998DiffeomorphismsGA,Bruveris2009TheMM} are examples of directly solving \eqref{eq: global optimisation} to obtain a diffeomorphism $T$ belonging to the same connected component as the identity map in the group of diffeomorphisms $\mathcal{D}(\Omega)$. Since an external forcing is not considered, these methods can be stated to be simplified versions of that proposed by \cite{Trouv2005MetamorphosesTL}. The regularization term in \eqref{eq: global optimisation} is chosen to be  
\begin{align}
\|T\|^2 = \int_{0}^1\|v(t)\|_V^2 dt,
\end{align}
for some prescribed norm $V$ and $T\approx \bfPhi_{\bfu_{n-1}}(\epsilon)\circ \cdots\circ \bfPhi_{\bfu_{0}}(\epsilon)$ as described in the introduction. Note that the optimisation problem in \cite{Beg2005ComputingLD,Bruveris2009TheMM} is a specific case of Eq.(6) in \cite{Trouv1995} and Eq.(1) in \cite{Trouv1998DiffeomorphismsGA}. 
%While methods in \cite{Beg2005ComputingLD,Bruveris2009TheMM} are more computationally demanding than the proposed method, the solution of the proposed method has a larger chance to be trapped in local minimums. 
It must be noted that the general theory of \cite{Trouv1995} does not particularly assume that $T^\#S$ is simply the composition, and that the sub-optimal solution in section 4.2 of \cite{Trouv1995} is closely related to the optimisation problem defined in \eqref{eq: optical flow for tensor fields, compact no boundary} and \eqref{eq: optical flow for tensor fields, compact with boundary}. Despite the fact that $T^\#$ in \cite{Beg2005ComputingLD,Bruveris2009TheMM,Trouv1995,Trouv1998DiffeomorphismsGA, Trouv2005MetamorphosesTL} is not particularly designed for geophysical fields and sometimes the boundary condition is ignored, the framework proposed in this manuscript can be thought of as a specific case of the general theory developed in \cite{Trouv1995}. But again, the main point of this paper is to show the importance of $T^\#$ when applying field alignement algorithms to geophysical fields. 

\subsection{About the boundary condition of \eqref{eq: optical flow for tensor fields, compact with boundary}}
Boundary condition is necessary in the optimisation problem \eqref{eq: optical flow for tensor fields, compact with boundary}. In fact, the boundary condition can be obtained by solving for \eqref{eq: optical flow for tensor fields, compact no boundary} or \eqref{eq: optical flow for tensor fields, compact with boundary} based on the data on $\partial\Omega$. The process can be illustrated using the following idealized example. 

Suppose that $\Omega$ is a three dimensional ocean, the boundary of which consists of two parts $\Omega_b$ and $\Omega_s$, in which $\Omega_b$ refers to the ocean basin (i.e. the land boundaries), and $\Omega_s$ the sea surface. To determine the boundary condition of $\bfu$ on $\Omega_b$ and $\Omega_s$, it is first natural to set $\bfu\big{|}_{\Omega_b} = 0$. %However, it is not reasonable to assume the same for $\Omega_s$. 
On $\Omega_s$, we must first solve for \eqref{eq: optical flow for tensor fields, compact with boundary} where $\theta_1$ and $\theta_2$ are tensor fields on $\Omega_s$. Boundaries of $\Omega_s$ coincide with the coast lines, and natural to set boundary conditions to be $0$. The domain $\Omega_s$ is a sub-manifold of the sphere. Hence the "Riemannian" context is necessary for this example. The complete process at each iterative step is illustrated in Algorithm 1.%\ref{alg: how to determine the bc}.

\begin{figure}[h!]
	\centering
	\includegraphics[width=0.9\linewidth]{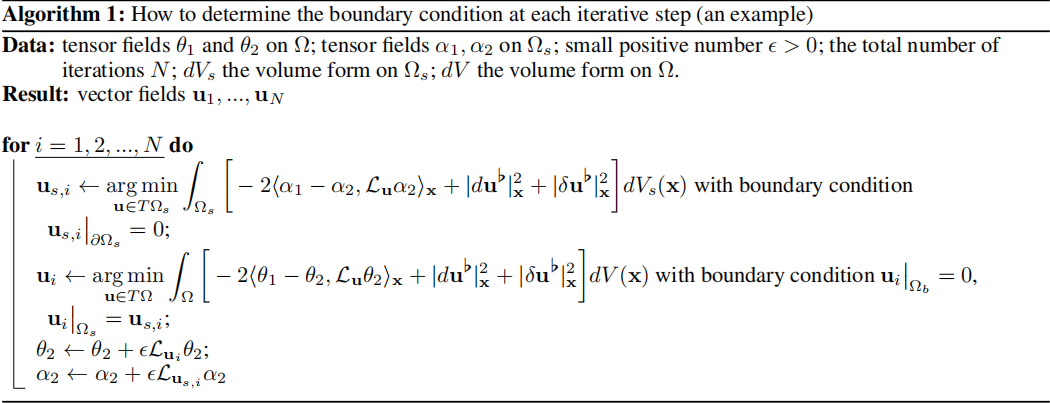}
	\label{alg: how to determine the bc}
\end{figure}

\subsection{The case for localized observations}

In practice, it is a common situation that some state variables are solely observed in a subdomain $\Omega_1\subset\Omega$, instead of in the full domain $\Omega$ (i.e. over a satellite swath). A weight function $W$ can be constructed and the first terms in Eq.\eqref{eq: optical flow for tensor fields, compact no boundary} and \eqref{eq: optical flow for tensor fields, compact with boundary} replaced by:
\begin{align}
\frac{\partial |\theta_1-\bfPhi_\bfu(s)_*\theta_2|_\bfx^2}{\partial s}\big{|}_{s=0}W(\bfx),
\end{align}
in which the pre-chosen $W(\bfx)$ takes value $1$ for a majority of points inside  $\Omega_1$ but decreases to $0$ smoothly as $\bfx$ approaches the boundary of $\Omega_1$. Then Eq.\eqref{eq: optical flow for tensor fields, compact no boundary, v2} changes to 
\begin{align}
\bfu = \argmin_{\bfu\in T\Omega} \int_{\Omega}\bigg{[}2\inp{\theta_1-\theta_2}{\mathcal{L}_\bfu\theta_2}_\bfx W(\bfx) + a_1|d\uflat|_\bfx^2 + a_1|\delta \uflat|_\bfx^2 + a_0|\uflat|_\bfx^2\bigg{]} dV(\bfx),\label{eq: optical flow for tensor fields, compact no boundary, local obs}
\end{align}
and Eq.\eqref{eq: optical flow for tensor fields, compact with boundary, v2} changes to 
\begin{align}
\begin{cases}
\bfu = \displaystyle\argmin_{\bfu\in T\Omega} \int_{\Omega}\bigg{[}2\inp{\theta_1-\theta_2}{\mathcal{L}_\bfu\theta_2}_\bfx W(\bfx) + a_1|d\uflat|_\bfx^2 + a_1|\delta \uflat|_\bfx^2 \bigg{]} dV(\bfx)\\
\bfu\big{|}_{\partial\Omega} = i_*\bfv.
\end{cases}\label{eq: optical flow for tensor fields, compact with boundary, local obs}
\end{align}
The theorem of existence and uniqueness of solution still holds for minimisation problems \eqref{eq: optical flow for tensor fields, compact no boundary, local obs} and \eqref{eq: optical flow for tensor fields, compact with boundary, local obs} following similar arguments given in the appendix.

\subsection{A displacement-based formulation of the nudging method in continuous time and continuous space}\label{subsection: displacement-based nudging}
With the notations and techniques introduced in section 2, a continuous-time nudging method the $K$ for which is based on the displacement of state variables, which has been mentioned in the introduction, can be explicitly formulated as the following:
\begin{align}
\frac{\partial \theta_x}{\partial t} = M_{\theta_x}(t, \theta_x) + \mathcal{L}_{\bfu(t)}\theta_x, \label{eq: displacement based nudging}
\end{align}
in which $\theta_x$ is the tensor field associated with the state variable $x$, and $M_{\theta_x}$ the corresponding forecast model for $\theta_x$, $\mathcal{L}_{\bfu(t)}\theta_x$ the Lie derivative of $\theta_x$ with respect to $\bfu(t)$, and $\bfu(t)$ is determined by the solution of optimisation problem \eqref{eq: optical flow for tensor fields, compact no boundary} or \eqref{eq: optical flow for tensor fields, compact with boundary}.   For penalty terms used in the optimisation problems \eqref{eq: optical flow for tensor fields, compact no boundary} and \eqref{eq: optical flow for tensor fields, compact with boundary}, it is explained in appendix \ref{appendix: the linearity of L_u theta on y-xb} that $\mathcal{L}_{\bfu(t)}\theta_x$ is indeed linear with respect to $\theta_1-\theta_2$, hence there exists an operator $\tilde{K}$ that acts linearly on $\theta_1-\theta_2$ such that $\tilde{K}(\theta_1 - \theta_2) = \mathcal{L}_{\bfu(t)}\theta_x$. $\tilde{K}$ depends on $x$ and $y^o-h(x^b)$. The discretised version of $\tilde{K}$ would be the matrix $K$ used in Eq.\eqref{eq: nudging in general}. Hence Eq.\eqref{eq: displacement based nudging} is a particular case of Eq.\eqref{eq: nudging in general}. Note that $\bfu(t)$ is estimated using the observed state variables $y^o(t)$. Then the same displacement vector field is applied to transport the unobserved state variables. Thus physically-consistent estimation of $\bfu(t)$ is of particular importance.

\section{Instantaneous nudging based on displacement error and its application in data assimilation}
\subsection{Methodology}
Let $S$ be the full state variable, $Y = h(S)$ some state variable derived from $S$. Assume that we have the model estimate $S_{\text{model}}$, and that $Y$ is fully observed on the domain. Our target is to 
\begin{itemize}
	\item[] (1) derive a displacement flow $\bfPhi$ so that $\bfPhi^\#Y^{\text{model}}\approx Y^{\text{obs}}$;
	\item[] (2) apply $\bfPhi$ to $S^{\text{model}}$ to correct the displacement of the full state variable $S$.
\end{itemize}
From section 2, we need to separately find the tensor fields for $Y$ and $S$. Consistent definitions should be determined by the dynamical equation of the system, leading to the determination of $T^\#$, hence the explicit formulation of Eq.\eqref{eq: optical flow for tensor fields, compact no boundary, v2} (or \eqref{eq: optical flow for tensor fields, compact with boundary, v2}).

Already indicated in \cite{Beezley2007MorphingEK} and \cite{Ravela2007DataAB}, this nudging strategy can be incorporated with ensemble Kalman filter. The basic idea is to correct the displacement of each ensemble member before applying EnKF. Methods reported in \cite{Beezley2007MorphingEK} and \cite{Ravela2007DataAB} use different cost functions to estimate the displacements. Both methods implicitly assume $T^\#S = S\circ T$ or $S\circ T^{-1}$. However, $T^\#$ determines how state variables are transported along the displacement vector fields. To avoid destroying the physical balance between state variables, only a proper definition of $T^\#$, for each state variable, can  ensure a correct transport along the same displacement vector field.

Simply illustrated, Fig.\ref{fig: 2}, considering $T^\#$ without prescribing the dynamics could possibly destroy the intrinsic feature of the unobserved physical fields. 
%It is assumed here that $T$ is obtained without using any information of the wind field. We remark that 
The definition of $T^\#$ is thus a fundamental component in all methods/algorithms which involve displacement of physical fields. In the Eq. (17) in \cite{Ravela2007DataAB}, the term $X^f_s({\bf p}_s)$ represents the $X^f_s$ adjusted by the displacement map ${\bf p}_s$. It can be equivalently written as $X^f_s\circ T$, for $T = {\bf p}_s$. For such a case, it can be suggested to replace $X^f_s({\bf p}_s)$ by ${\bf p}_s^\#X^f_s$, and similarly for all terms wherever $T^\#S$ is implied. The specific definition of ${\bf p}_s^\#$ should then depend on the state variables in $X^f_s$. Since the calculation of $T^\#$ is straightforward, replacing $S\circ T$ by $T^\#S$ would not dramatically increase the computational cost of the original algorithm. Note, there are already several attempts to apply this alignment method \cite{Ravela2007DataAB} in more practical observing system simulations experiments (OSSEs) \cite{Ravela2014SpatialIF, Jankov2021PartitionOF}. In all these experiments, the fields of all concerned state variables, although noisy and sometimes sparse, are used as the reference. 
%However, for online applications it can not be assumed that all state variables are sparsely observed. 
Yet, even observing all state variables, the alignment method \cite{Ravela2007DataAB} may not always provide the optimal solution. 
%If only a few state variables are observed, this method would then possibly produce physically imbalanced results. 
Detailing a numerical demonstration to improve this method with a physically-driven $T^\#$, will be the focus of a future paper. 
%Thus, it questions the application of the displacement map $T$ estimated in \cite{Beezley2007MorphingEK} and \cite{Ravela2007DataAB} to modify the displacement of unobserved physical fields.

To demonstrate the advantage of physically consistent $T^\#$, the following simple version of morphed EnKF is numerically compared with the plain EnKF algorithm:
\begin{itemize}
	\item[] (1) choose tensor fields for the observed state variable $Y$ and the full state variable $X$;
	\item[] (2) find the displacement flow $\bfPhi_i$ according to the observation $Y^{\text{obs}}$ and the model estimates $Y^{\text{model}}_i$ from the $i-$th ensemble member;
	\item[] (3) apply $\bfPhi_i$ to the full state vector of the $i-$th ensemble member: $X^f_i \gets \bfPhi_i^\# X^f_i$;
	\item[] (4) apply plain EnKF to the updated ensemble $\{X^f_i\: i=1,...,N_e\}$.
\end{itemize}
A morphed ensemble member based on the plain $T^\# S = S\circ T^{-1}$ is also presented for comparison.

\subsection{Numerical results}
The data assimilation experiment is conducted for only one time step, using the thermal shallow water equation \cite{Warneford2013TheQT}. This model consists of three state variables: $h-$the water height, $\bfv = (v_1,v_2)-$the velocity field, and $\Theta-$the buoyancy (or density contrast):
\begin{align}
\frac{\partial h}{\partial t} + \nabla\cdot(h \bfv) &= 0,\label{eq: thermal shallow water, h}\\
\frac{\partial \Theta}{\partial t} + (\bfv\cdot\nabla)\Theta &= -\kappa(h\Theta - h_0\Theta_0),\label{eq: thermal shallow water, Theta}\\
\frac{\partial \bfv}{\partial t} + (\bfv\cdot\nabla)\bfv + f\hat{\bf z}\times \bfv &= -\nabla(h\Theta) + \frac{1}{2}h\nabla\Theta.\label{eq: thermal shallow water, v}
\end{align}
Both $\Theta$ and $h$ are assumed strictly positive at each point. It is assumed that both the absolute vorticity $\omega = \frac{\partial v_2}{\partial x} - \frac{\partial v_1}{\partial y}$ and $h$ are fully observed, while $\Theta$ is completely unobserved. The details of the experimental setup can be found in the appendix. A Python code to completely reproduce the numerical results is provided at http://doi.org/10.5281/zenodo.10252175.

\begin{figure}
	\centering
	\includegraphics[width=1.0\textwidth]{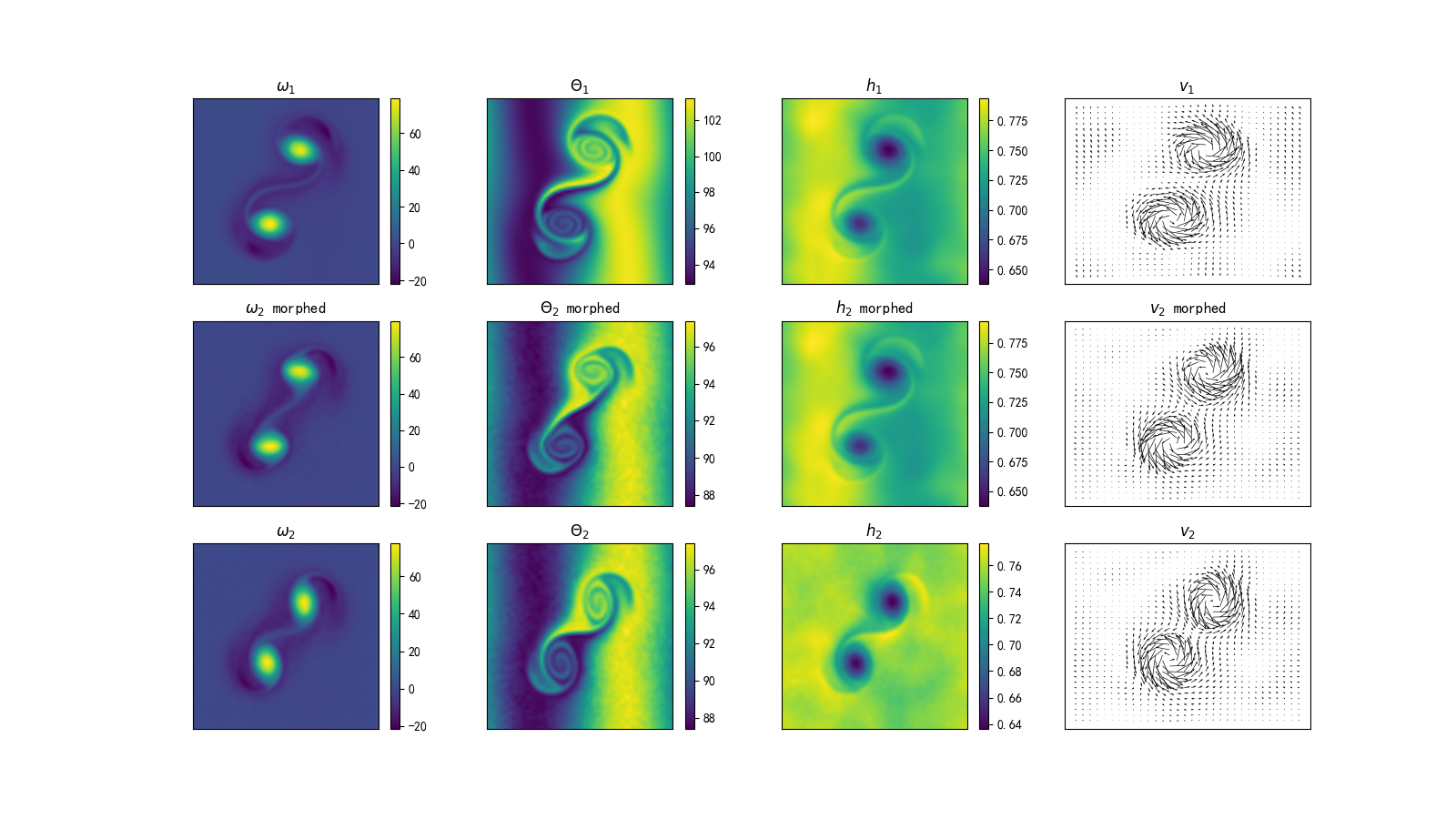}
	\caption{Comparison of the target fields (the first row), the original fields (the third row), and the morphed fields (the second row). This is the result of one of the ensemble members. The displacement map is calculated based on the $\omega$-field and the $h$-field. Then the displacement map is applied to the buoyancy field $\Theta$ and the velocity field $\bfv$. For this ensemble member, it is not hard to see that both the buoyancy and the velocity are partially aligned with the truth. However, the results for another ensemble member is not so optimistic. See Fig.\ref{fig: xfmember2}.  }
	\label{fig: xfmember1}
\end{figure}

\begin{figure}
	\centering
	\includegraphics[width=1.0\textwidth]{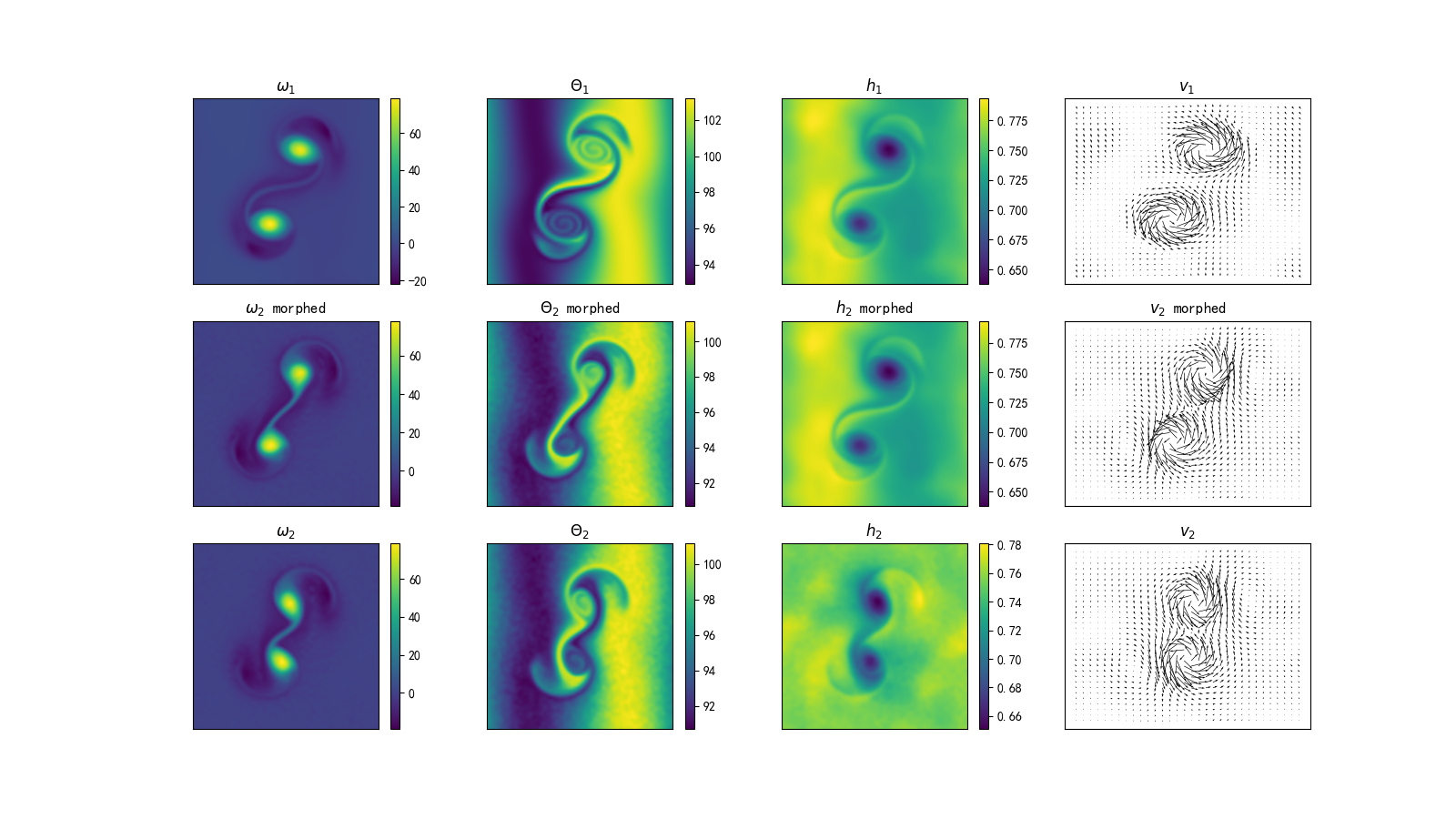}
	\caption{ Comparison of the target fields (the first row), the original fields (the third row), and the morphed fields (the second row). This is the result for another ensemble member. The displacement map is calculated based on the $\omega$-field and the $h$-field. Then the displacement map is applied to the buoyancy field $\Theta$ and the velocity field $\bfv$. For this ensemble member, the initial discrepancy between $\Theta_2$ and $\Theta_1$ is larger than that in Fig.\ref{fig: xfmember1}. And the result is not as good as the result for the ensemble member shown in Fig.\ref{fig: xfmember1}.  }
	\label{fig: xfmember2}
\end{figure}
Fig.\ref{fig: xfmember1} and Fig.\ref{fig: xfmember2} show two examples of the (interpolated) truth (the target fields), the prior estimate of one member (the original fields), and the morphed prior estimate (the morphed fields). As expected, the phase of the original vortex can be adjusted, to some extent, but not perfectly. This is not surprising because the cost function in the optimisation problem \eqref{eq: optical flow for tensor fields, compact no boundary} (or \eqref{eq: optical flow for tensor fields, compact with boundary}) has two parts. The first part is derived from the dynamics of the original system. But the second part is  a regularization term which usually is not based on the original dynamics of the system. Thus, the virtual displacement flow $\Phi$ cannot transform $\Theta_2$ to exactly match $\Theta_1$. However, the correct choice of tensor fields maintains, to some extent, the dynamical balance of the three fields during the morphing process.
%to some extent. 
Looking more closely, Fig.\ref{fig: xfmember2}, the difference between the morphed $h_2$ field and the target $h_1$ field is much smaller than between the morphed $\omega_2$ and the target $\omega_1$. Numerically, it is also found that the $h$-field is easier to align than the $\omega$-field. It needs further exploration if there is any particular reason for this.
%We don't have a theoretical explanation of this fact. 

\begin{figure}
	\centering
	\includegraphics[width=1.\textwidth]{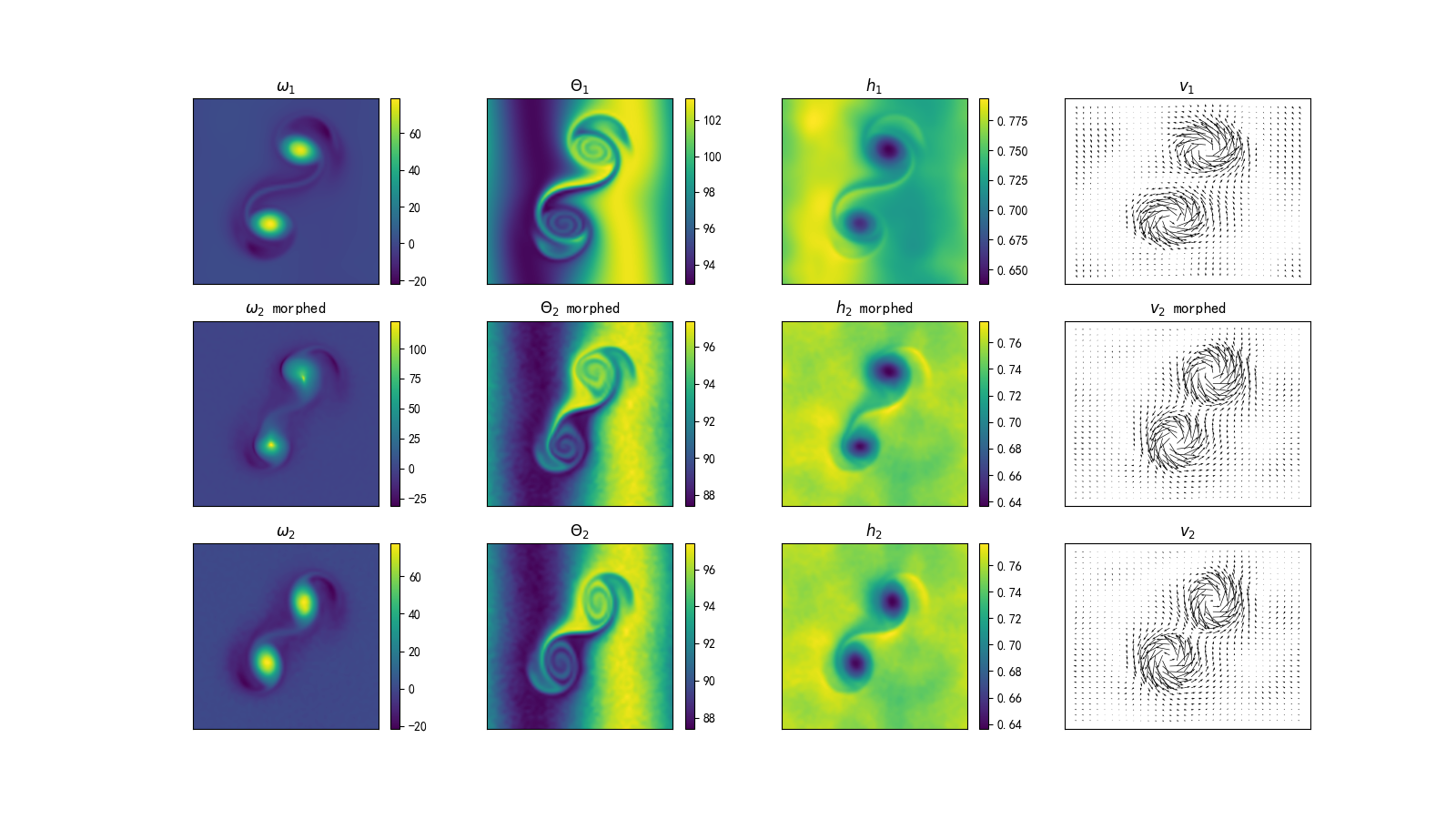}
	\caption{The target field (first row), the original field (third row), and the field morphed by $T^\#S = S\circ T^{-1}$  for $S = h,\omega, \Theta, \bfv$ (second row). Since $T^{\#}S = S\circ T$ is not physically consistent, a singularity has been generated in the morphed $\omega_2$-field.}
	\label{fig: morphed with no DG}
\end{figure}

To demonstrate the usefulness of $T^\#$, the nudging process, without introducing the concept of tensor fields, is also considered and tested. In this case, $T^\#S = S\circ T^{-1}$, independent of $S = h$, $\omega$, $\Theta$, or $\bfv$. Initial values, $\omega_2$ $h_2$, etc, correspond to fields of one of the ensemble priors. The morphing process is first run for $N=10000$ time steps. For this specific member, the mean-squared-error (MSE) of $h$ always decreases, while the MSE for $\omega$ decreases till $N\approx 6000$ and then starts to increase. Fields at $N=6000$ are plotted Fig.\ref{fig: morphed with no DG}.  It can be seen that a singularity has been generated in the morphed $\omega_2$-field.

\begin{figure}
	\centering
	\includegraphics[width=1.\textwidth]{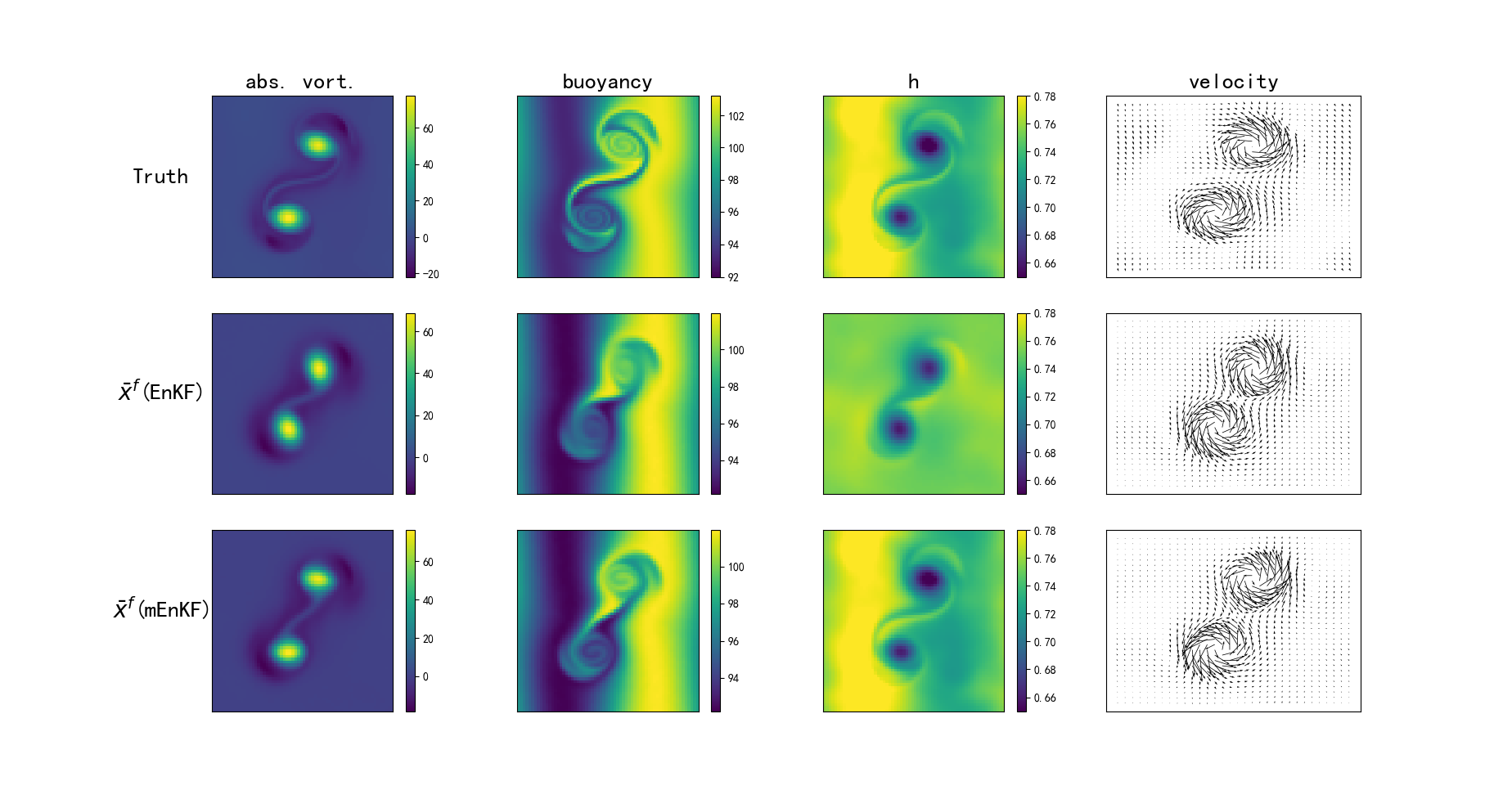}
	\caption{The truth (first row), the prior mean of ensemble (the second row), and the prior mean of the morphed ensemble (the third row). Some small scale features of the vortex is lost in the arithmetic mean of the original ensemble members. However, the arithmetic mean of the morphed ensemble members still preserve some of the small scale features.} 
	\label{fig: xfmean}
\end{figure}
Fig. \ref{fig: xfmean} presents the truth, the prior mean of the ensemble, and the prior mean of the morphed ensemble. The direct prior mean has completely lost the small scale structures inside the vortex, while the prior mean of the morphed ensemble still maintain, to some extent, the small scale features. Stated in the introduction, $y^o - x^b$ may no longer be a good representation of the error. Instead, the location error, or more generally the displacement flow, can better represent the error of each ensemble member.  Intuitively, the state vectors represented by the ensemble members are supposed to lie on a curved manifold rather than on a flat Euclidean space. In this case, it is not surprising to see that the arithmetic mean lies outside of the manifold, which results in structural error of the whole field. To address this problem, the Fr\'{e}chet-type of mean instead of the arithmetic mean should be used to define the ensemble mean.

\begin{figure}
	\centering
	\includegraphics[width=1.\textwidth]{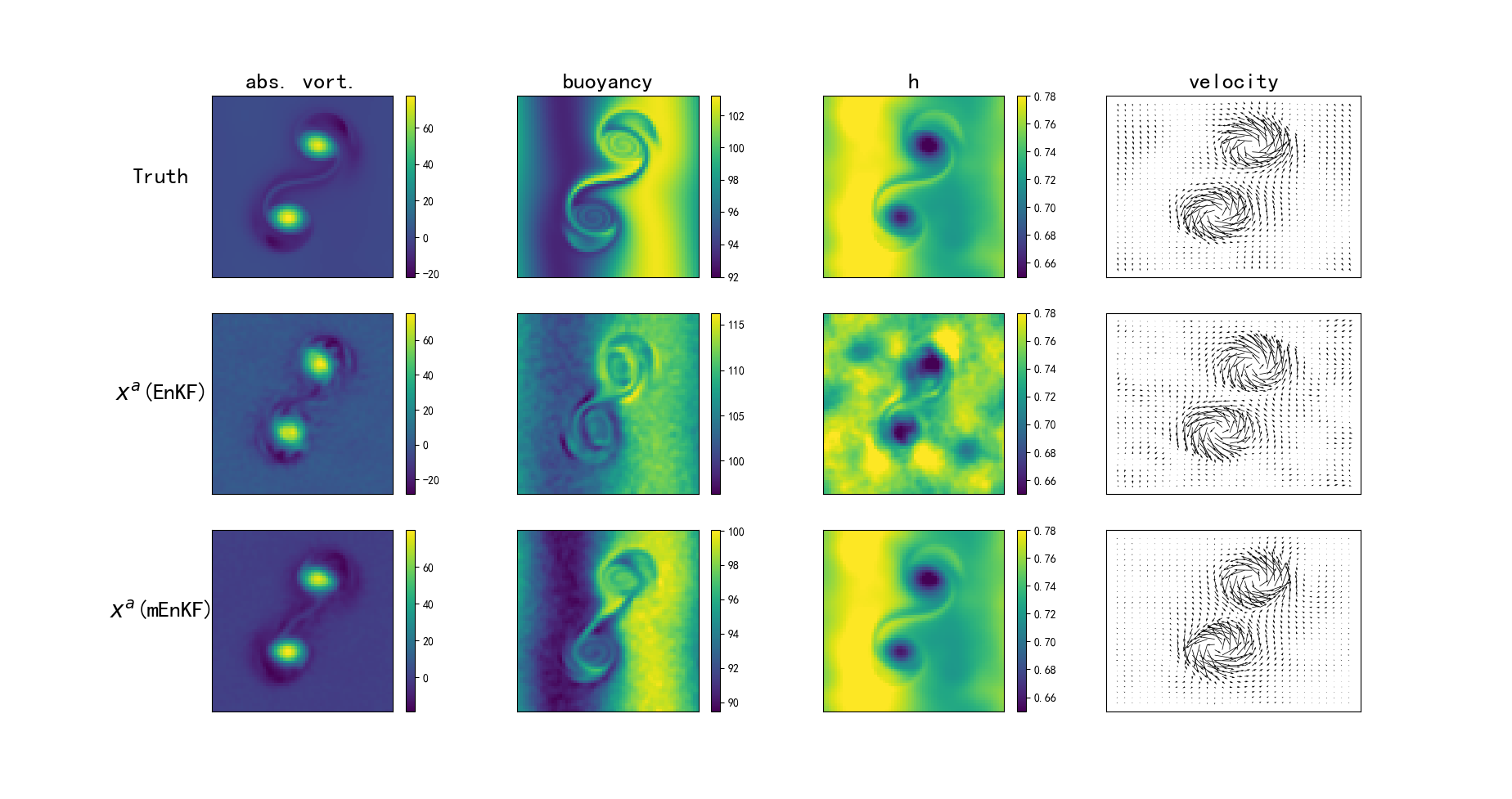}
	\caption{The truth (first row), the posterior estimate of one member (the second row), and the posterior estimate of the corresponding morphed member (the third row). The plain EnKF algorithm largely destroys the intrinsic feature of the double vortex field, while the morphed EnKF method could still preserve some of the intrinsic feature, although not perfectly. }
	\label{fig: xamember}
\end{figure}
Fig.\ref{fig: xamember} presents the posterior estimate of one ensemble member using EnKF and morphed EnKF (denoted by mEnKF for short). Apparently, the plain EnKF results in erroneous estimates, while the mEnKF still produces reasonable estimates with fine scale features inside the vortex. Again, this is the disadvantage of linear algorithms when the ensemble members do not lie on a linear subspace of the Euclidean space of state vectors.

\section{Conclusion}
It is demonstrated that physically consistent displacements must be considered when field alignment algorithms are applied to geophysical fields. How a physical field morphs along a vector field must then be specifically designed according to the underlying physical principles. For such a purpose, differential geometry appears to be a proper tool to mathematically formulate the displacement flow between two physical fields. It is mathematically proved that the solution of the optimisation problem always exists and is unique. When a time sequence of snapshots of geophysical fluids is available, a reformulation of the generalized OF algorithm \cite{Wildes2000RecoveringEO} can also be obtained from this differential geometry perspective.

Despite that boundary condition is not considered in \cite{Trouv1995}, the proposed framework appears as a special case of the general theory developed in \cite{Trouv1995}. However, \cite{Trouv1995} and its subsequent applications/methods \cite{Beg2005ComputingLD,Bruveris2009TheMM,Trouv1998DiffeomorphismsGA,Trouv2005MetamorphosesTL} never explicitly defined $T^\#$ for geophysical state variables. Correspondingly, relevant methods in the geophysical community \cite{Ravela2007DataAB,Beezley2007MorphingEK} seemingly failed to endorse the importance of $T^\#$. Here demonstrated, $T^\#$ is explicitly defined, with the help of differential geometry, for certain state variables. A key advantage is to consistently correct the displacement of some state variable while best maintaining the intrinsic structure of the underlying physical fields. Methods in \cite{Ravela2007DataAB,Beezley2007MorphingEK} can thus likely be improved by incorporating the correct $T^\#$. 
%We reserve the thorough study of this aspect to the future.

The proposed framework can also be used in a nudging process when only part of the state variables are observed. It can also be used to correct the displacement of each ensemble member before applying the linear EnKF algorithm. Numerical results with a double-vortex model show that the morphed EnKF algorithm produces more reliable posterior estimates than the plain EnKF algorithm. The preferable usage of $T^\#$ is well demonstrated by comparing fields morphed by algorithms derived or not from a more careful definition of $T^\#$.

It must be pointed out that all physical fields could not always be associated to a tensor field. The use of a tensor field is only required to define the correct $T^\#$. For an arbitrary derived physical field $S$, $T^\#S$ should be determined by the forward operator $S = F(X)$, where $X$ is the full state vector. Thus, for a specific derived variable $S$, whether a displacement vector field can be determined by the optimisation procedure \eqref{eq: optical flow for tensor fields, compact no boundary} or \eqref{eq: optical flow for tensor fields, compact with boundary} needs to be further studied.

\appendix

\section{Differential forms, pull-back map $T^*$, Hodge star operator and Lie derivatives in Euclidean spaces for practitioners in data assimilation}\label{Appendix: differential geometry concepts}

For a thorough introduction of these concepts in general Riemannian manifolds, please refer to \cite{Chern1999LecturesDiffGeo} or any other textbook in differential geometry.  Suppose $\Omega$ is a $n$-dimensional manifold in $\mathbb{R}^n$. Intuitively this means that $\Omega$ is an open subset of $\mathbb{R}^n$, instead of a curved surface in $\mathbb{R}^n$. For instance, ignoring the general relativity effect, the whole atmosphere is such a manifold. But the Earth/sea surface or the boundary of atmosphere is not such a manifold. The concept of "Riemannian" in this case is hidden behind the Euclidean assumption. Hence in this appendix we will not explicitly point out how the Riemannian metric induces the inner product on certain spaces. Because everything just follows the natural sense. 

\subsection{Vector field}
For each point $\bfx\in\Omega$, the tangent space (or tangent plane) at $\bfx$ is denoted by $T_\bfx\Omega$. $T_\bfx\Omega$ is another $\mathbb{R}^n$ which is associated to the point $\bfx$ only. Note that the tangent space for different points $\bfx_1$ and $\bfx_2$ are different, although both of them have the same structure as $\mathbb{R}^n$.  Although $T_{\bfx_1}\Omega$ and $T_{\bfx_2}\Omega$ are different for different points $\bfx_1$ and $\bfx_2$, vectors in $T_{\bfx_1}\Omega$ and $T_{\bfx_2}\Omega$ can be 1-1 represented using the same coordinate system, i.e the coordinate system for the $\mathbb{R}^n$ that $\Omega$ lives in. However, this is not the true for general Riemannian manifolds.  Thus a vector field $\bfv$ on $\Omega$ can be expressed as $\bfv({\bfx}) = \sum_{i=1}^na_i(\bfx)\bfe_i$. 

We only care about vector fields $\bfv$ that have some regularity, for instance $H^0$ (the components of $\bfv$ are square-integrable),  $C^0$ (continuous), $C^1$ ( all $a_i(\bfx)$ are continuously differentiable), smooth (all $a_i$'s are infinitely many times differentiable), or $H^1$ (differentiable, $L^2$, and the partial derivatives of $a_i$ are $H^0$). Note that $H^1$ and $C^1$ are not equivalent. The set of vector fields that enjoy the regularity of $H^0$ or $H^1$ are denoted by $H^0(T\Omega)$ and $H^1(T\Omega)$, respectively. 

\subsection{Differential $1$-form}
Differential $1$-forms are dual elements of vector fields. Suppose $\omega$ is a differential $1$-form. This means that  for each point $\bfx\in\Omega$, $\omega(\bfx)$ is a linear function on $T_{\bfx}\Omega$. There are a lot of linear functions on $T_\bfx\Omega$. For instance, take any vector $\bfu\in T_\bfx\Omega$, the natural inner product $\inp{\bfu}{\cdot}$ is a linear function on $T_\bfx\Omega$. In fact, all linear functions on $T_\bfx\Omega$ can be represented in such a way. This means that each linear function $l(\cdot)$ on $T_{\bfx}\Omega$ corresponds to a vector $\bfu_l\in T_\bfx\Omega$. Thus the space of all linear functions on $T_{\bfx}\Omega$, denoted by $A$ just for now, seems can be identified with $T_\bfx\Omega$ itself. However, we do not identify $A$ with $T_\bfx\Omega$ even in the Euclidean case. The reason is explained in the introduction of pull-back maps. The space of all linear functions on $T_\bfx\Omega$ is called the dual space of $T_\bfx\Omega$, denoted by $T_\bfx^*\Omega$. Again, $T^*_\bfx\Omega$ and $T_\bfx\Omega$ both look like $\mathbb{R}^n$. But they are different. $T_\bfx^*\Omega$ is called the cotangent space at $\bfx$.

Let $\{\bfe_i\}_{i=1,...,n}$ be the natural orthogonal basis of the space $\mathbb{R}^n$ that $\Omega$ lives in. This can also be viewed as a basis for $T_\bfx\Omega$ for any $\bfx\in\Omega$. However, it can not be viewed as a basis for $T^*_\bfx\Omega$. We use $\{\bfe_i^*\}_{i=1,...,n}$ to denote the basis for $(\mathbb{R}^n)^*$, so that $\{\bfe_i^*\}_{i=1,...,n}$ is the dual basis of $\{\bfe_i\}_{i=1,...,n}$ in the following sense: $\bfe_i^*(\bfe_j) = 1$ if $i=j$ and $0$ otherwise. Similar to the tangent space, the cotangent space at different points share the same basis in the Euclidean case. Thus a differential $1$-form $\omega$ on $\Omega$ has the expression: $\omega(\bfx) = \sum_{i=1}^nb_i(\bfx)\bfe_i^*$, where $b_i(\cdot)$'s are functions on $\Omega$.

However, $\omega = \sum_{i=1}^nb_i\bfe_i^*$ is not the usual way of writing the differential $1$-forms. Usually we write it as $\omega = \sum_{i=1}^n b_idx_i$. The $dx_i$ here uses exactly the same symbol as the $dx^i$ used in calculus. But it is not the infinitesimal element $dx_i$ in Calculus. $dx_i$ here refers to $\bfe_i^*$. There are several reasons for  writing it in this way. $dx_i$ can be interpreted as the exterior derivative (to be defined) of the function $x_i$. Hence it is directly related to Calculus. For a differential $n$-form  $\omega =f(\bfx)dx_1\wedge\cdots\wedge dx_n$ (or $\omega = f(\bfx)\bfe_1^*\wedge\cdots\wedge\bfe_n^*$), it has a well-defined integral on the whole manifold: $\int_{\Omega}\omega := \int_{\Omega} f(\bfx)dx_1dx_2\cdots dx_n$, the right hand side of which is the integral in the sense of calculus. The wedge symbol $\wedge$ is defined for tensor product of linear spaces which is not directly relevant to calculus. Thus $dx_i$ is also used to represent $\bfe_i^*$.  

Similar to the vector fields, we only care about differential $1$-forms with certain regularity. The corresponding spaces are denoted by $H^0(T^*\Omega)$, $H^1(T^*\Omega)$, respectively.

\subsection{Differential $0$-form}
Differential $0$-forms on $\Omega$ are just functions on $\Omega$.
\subsection{Differential $k$-form}
A differential $k$-form $\omega$ can be defined based on the tensor product of cotangent spaces \cite{Chern1999LecturesDiffGeo}. In Euclidean space a differential $k$-form has the following expression: \begin{align}
\omega(\bfx) = \displaystyle\sum_{1\leq i_1<i_2<\cdots<i_k\leq n} f_{i_1i_2...i_k}(\bfx)dx_{i_1}\wedge\cdots\wedge dx_{i_n} \label{eq: Appendix: w expression in ascending index}
\end{align}
The wedge product $\wedge$ of a differential $k$-form and a differential $l$-form gives a differential $(k+l)$-form. There are some rules for the wedge product. These rules are not artificially made but are consequences of the definition. Here we do not provide the complete definition of differential $k$-forms. We only list the rules of  wedge product $\wedge$. To obtain a thorough understanding of why these rules are true, please refer to \cite{Chern1999LecturesDiffGeo}  or any other textbook on differential geometry.  The space of differential $k$-forms with $H^s$ regularity is denoted by $H^s(\Lambda^kT^*\Omega)$. Using the expression \eqref{eq: Appendix: w expression in ascending index}, in the Euclidean case the regularity of $\omega$ is simply equivalent to the regularity of $f_{i_1i_2...i_k}$. 

Assume that $\alpha$ is a differential $k$-form and $\beta$ a  differential $l$-form, where $ 0\leq k,l$ are integers. And $\gamma$ is any differential form. Then the following rule of wedge product holds:
\begin{align}
\alpha\wedge \beta = (-1)^{kl}\beta\wedge\alpha \\
(\alpha\wedge\beta)\wedge\gamma = \alpha\wedge (\beta\wedge\gamma) =: \alpha\wedge\beta\wedge\gamma
\end{align}
Note that $dx_i$ itself is a differential $1$-form. The above rule then implies that $dx_i\wedge dx_j = -dx_j\wedge dx_i$. This immediately implies that $dx_i\wedge dx_i = 0$. Thus if the subindex of $\alpha_I = dx_{i_1}\wedge\cdots\wedge dx_{i_k}$ and the subindex of $\alpha_J = dx_{j_1}\wedge\cdots\wedge dx_{j_k}$ are the same set but only differ by a permutation, then essentially $\alpha_I = \alpha_J$ or $\alpha_I = -\alpha_J$. This is why we can assume ascending subindex in the expression of $\omega$ in \eqref{eq: Appendix: w expression in ascending index}. Also, this rule implies that there can not be more than $n$ subindex in one term because those terms would automatically equal to 0. Thus all differential $(n+k)$-forms are automatically 0 for $k\geq 1$.  

\subsection{The integral of differential $n$-form}
If $\Omega$ is a $n$-dimensional Euclidean manifold, and $\omega$ is a differential $n$-form. Assume that $\omega(\bfx) = f(\bfx)dx_1\wedge\cdots\wedge dx_n$. Then the integral of $\omega$ on $\Omega$ is defined to be $\displaystyle\int_{\Omega}\omega := \displaystyle\int_{\Omega}f(\bfx)dx_1dx_2\cdots dx_n$, where the integral on the right is the integral in the sense of calculus. 

\subsection{The volume form $dV$}
For a $n$-dimensional manifold $\Omega$, the volume form, denoted by $dV$, is a differential $n$-form. In Euclidean case, $dV = dx_1\wedge\cdots\wedge dx_n$. It is called a volume form because for any $n$-dimensional sub-manifold $M$ (you can think of it as an open subset of $\Omega$), $\displaystyle\int_{M}dV$ is just the volume of $M$.

\subsection{The inner product of differential $k$-forms}
Let $\alpha,\beta$ be differential $k$-forms. $\alpha = \displaystyle\sum_{|I|=k}f_Idx_I$ and $\beta = \sum_{|I|=k}g_Idx_I$, where $I$ refers to any ascending sequence of $k$ different subindex. 
The inner product of $\alpha$ and $\beta$ at each point $\bfx\in\Omega$ is defined by: $\inp{\alpha}{\beta}_\bfx = \displaystyle\sum_{|I|=k}f_Ig_I$.

\subsection{Exterior derivative of differential $k$-forms ($d\omega$)}
Assume that $\omega$ is a smooth differential $k-$form: $\omega = \displaystyle\sum_{|I|=k}f_Idx_I$. The operator $d$ is called the exterior derivative. It maps a differential $k$-form to a differential $(k+1)$-form.
\begin{align}
d\omega = \sum_{I}\sum_{i=1}^n\frac{\partial f_I}{\partial x_i}dx_i\wedge dx_I.
\end{align}
It is possible that for different subindex $I_1$ and $I_2$ and different $i_1$ and $i_2$, the term $dx_{i_1}\wedge dx_{I_1} = \pm dx_{i_2}\wedge dx_{I_2}$. Thus similar terms need to be combined. The exterior derivative has the following rules:
\begin{align}
&d(\alpha\wedge \beta) = (d\alpha)\wedge \beta + (-1)^k\alpha\wedge (d\beta)\\
&dd\alpha = 0\label{eq: Appendix: dd = 0}
\end{align}
where $\alpha$ is a smooth differential $k$-form and $\beta$ is a smooth differential $l$-form. $k$ and $l$ can be any nonnegative integers. Due to Eq.\eqref{eq: Appendix: dd = 0}, people write $dd = 0$ or $d^2 = 0$. Stoke's formula has a nice and clean presentation using differential forms:
\begin{align}
\int_{\Omega}d\omega = \int_{\partial\Omega}\omega,
\end{align}
where $\omega$ is a smooth differential $(n-1)$-form, and $\partial\Omega$ refers to the boundary of $\Omega$ with proper orientation. $\partial\Omega$ is a manifold of dimension $n-1$. In most of the time $\partial \Omega$ can not be a manifold with Euclidean metric. 

\subsection{Hodge star operator ($*$), codifferential operator ($\delta$), and Hodge Laplacian ($\Delta$)}
All these are linear operators. We only list how Hodge star operator acts on differential forms when $n=2$ or $3$. A complete definition of these operators can be found at page 20 of \cite{Schwarz1995HodgeD}. 

We start with Hodge star operator $*$. Hodge star operator maps a differential $k$-form to a differential $(n-k)$-form. Basically if $\omega = \displaystyle\sum_{i=1}^mf_{I_i}dx_{I_i}$, where each subindex $I_i$ is a $k$-index. Then $*\omega$ has the expression  $*\omega = \displaystyle\sum_{i=1}^mg_{J_i}dx_{J_i}$, where each subindex $J_i$ is a $(n-k)$-index such that the set union $J_i\cup I_i = \{1,2,...,n\}$. And $g_{J_i} = \pm f_{I_i}$, where the sign depends on the order of $I_i$ and $J_i$. We list some examples in the following which would be enough for the calculation of $*\omega$ for 2- or 3-dimensional Euclidean case.

For $n=2$, $*dx_1 = dx_2$, $*1 = dx_1\wedge dx_2$, $*dx_2 = -dx_1$, $*dx_1\wedge dx_2 = 1$. Note that the constant $1$ is a function (hence a differential $0$-form) on $\Omega$.

For $n=3$, $*1 = dV$, $*dx_1 = dx_2\wedge dx_3$, $*dx_2 = -dx_1\wedge dx_3$, $*dx_3 = dx_1\wedge dx_2$, $*(dx_1\wedge dx_2) = dx_3$, $*(dx_1\wedge dx_3) = -dx_2$, $*(dx_2\wedge dx_3) = dx_1$, $*dV = 1$. 

The codifferential operator $\delta$ maps a differential $k$-form to a differential $k-1$ form, $k=0,1,2,...,n$. $\delta$ maps a function to $0$. When $\Omega$ is compact and has no boundary (for instance the sphere), $\delta$ is the adjoint operator of the exterior derivative $d$. More specifically, for smooth differential $k$-form $\alpha$ and smooth differential $(k+1)$-form $\beta$, $\displaystyle\int_{\Omega}\inp{d\alpha}{\beta}_{\bfx} dV(\bfx) = \displaystyle\int_{\Omega}\inp{\alpha}{\delta\beta}_\bfx dV(\bfx)$. However, this equation does not hold in general when $\Omega$ has boundary.

When $\delta$ acts on differential $k$-form, $\delta := (-1)^k*d*$, where $*$ is the Hodge star operator and $d$ is the exterior derivative operator. Similarly, we have $\delta\delta = 0$. The Hodge Laplacian operator is defined to be $\Delta  := d\delta + \delta d$.

\subsection{The push-forward operator $T_*$}
Suppose that $\Omega_1$ and $\Omega_2$ are two $n$-dimensional (Euclidean) manifolds. A diffeomorphism $T:\Omega_1\to\Omega_2$ is a smooth map which is 1-1, onto, and whose inverse is also smooth. Since we are looking at Euclidean manifolds, $T$ has $n$ components: $T = (T_1,...,T_n)$. Each $T_i$ is a function on $\Omega_1$. Suppose $\bfp\in\Omega_1$ and $\bfq\in\Omega_2$ are two points such that $T(\bfp) = \bfq$. Then $T$ induces a map $T_*: T_\bfp\Omega_1\to T_\bfq\Omega_2$. Think of any smooth curve in $\Omega_1$, $\gamma = \{\gamma(t): 0<t<1\}$, that passes through $\bfp$. The image of $\gamma$ would be a curve that passes through $\bfq$. Denote by $\bfu$ the tangent vector of $\gamma$ at $\bfp$, and by $\bfv$ the tangent vector of $T\circ\gamma$ at $\bfq$. Then $T_*$ is defined to map $\bfu$ to $\bfv$. The map $T_*$ can be calculated explicitly using the Jacobian matrix of $T$:
\begin{align}
(T_*\bfe_i)\Big{|}_\bfq = \sum_{j=1}^n\frac{\partial T_j}{\partial x_i}\Big{|}_{\bfp}\bfe_j.
\end{align}
If $\bfu(\bfx) = u_1(\bfx)\bfe_1 + ... + u_n(\bfx)\bfe_n$ is a vector field on $\Omega_1$, $T_*$ can not map $\bfu$ to a vector field at $\Omega_2$ for general smooth map $T$. Suppose $T$ maps two points $\bfp_1$ and $\bfp_2$ to the same $\bfq$, then it would be unclear how to define  $T_*\bfu$ at $\bfq$. However, $T_*$ is well-defined if $T$ is a diffeomorphism. The formula is the following:
\begin{align}
(T_*\bfu)\Big{|}_{\bfq} = \sum_{i=1}^nu_1(T^{-1}(\bfq))(T_*\bfe_i)\Big{|}_{\bfq}
\end{align}

\subsection{The pull-back operator $T^*$}
For $\Omega_1$, $\Omega_2$ and $T$ defined above, the pull-back operator $T^*$ maps a cotangent vector $\beta$ at $\bfq$ to a cotangent vector $\alpha = T^*\beta$ at $\bfp$ so that $\alpha(\bfu) = \beta(T_*\bfu)$ for any $\bfu\in T_\bfp\Omega_1$. Assume that $\Omega_2$ has coordinates $ (y_1,...,y_n)$. Using the basis $\{dx_1,...,dx_n\}$ for $T^*_\bfp\Omega_1$ and $\{dy_1,dy_2,...,dy_n\}$ for $T^*_\bfq\Omega_2$, $T^*$ has the following expression:
\begin{align}
(T^*dy_i)\Big{|}_{\bfp} = \sum_{j=1}^n\frac{\partial T_i}{\partial x_j}\Big{|}_{\bfp}dx_j
\end{align}
Let $\beta$ be a differential $1$-form on $\Omega_2$, such that $\beta(\bfy) = \displaystyle\sum_{i=1}^n g_i(\bfy)dy_i$ for any $\bfy\in\Omega_2$. Then for any $\bfx\in\Omega_1$, 
\begin{align}
(T^*\beta)(\bfx) = \sum_{i=1}^ng_i(T(\bfx))\frac{\partial T_i}{\partial x_j}(\bfx)dx_j 
\end{align}
Note that $T^*$ is well-defined for any smooth map $T$. $T^*$ then generalizes to differential $k$-forms using the following rule:
\begin{align}
T^*(\beta_1\wedge\beta_2)) = (T^*\beta_1)\wedge (T^*\beta_2)
\end{align}
for any differential forms $\beta_1$ and $\beta_2$ on $\Omega_2$. More explicitly, 
\begin{align}
T^*(dy_{i_1}\wedge\cdots\wedge dy_{i_k}) = \sum_{|J| = k}\det\Big{[}\frac{\partial (T_{i_1},...,T_{i_k})}{\partial (x_{j_1},...,x_{j_k})}\Big{]}dx_{j_1}\wedge\cdots\wedge dx_{j_k}.
\end{align}
Another important rule for $T^*$ is that it commutes with the exterior derivative $d$:
\begin{align}
dT^* = T^* d.
\end{align}

\subsection{On the generalization of $T_*$ to differential forms}\label{Appendix: explain T_* acts on differential forms}
Suppose that $\Omega_1$ and $\Omega_2$ are two different Euclidean manifolds and $T:\Omega_1\to\Omega_2$ is a diffeomorphism. Let $\alpha$ be a differential $k$-form on $\Omega_1$ and $\bfu$ a vector field on $\Omega_1$. We want to generate a differential $k$-form $\beta$ on $\Omega_2$ and a vector field $\bfv$ on $\Omega_2$. Almost immediately we choose $\bfv = T_*\bfu$ and $\beta = (T^{-1})^*\alpha$. Now assume that $\Omega_1$ and $\Omega_2$ are the same manifold and $T:\Omega\to\Omega$ is a displacement map (a diffeomorphism). Then for a given vector field $\bfu$ and a given differential form $\omega$ on $\Omega$, following the same principle we should define the displaced $\bfu$ to be $T_*\bfu$ and the displaced $\omega$ to be $(T^{-1})^*\omega$. This is why we generalize $T_*$ in the way written in section 2 of the main text.

\subsection{A practical way to calculate the Lie derivative of a differential form}
Given a vector field $\bfu(\bfx) = u_1(\bfx)\bfe_1 + \cdots + u_n(\bfx)\bfe_n$ and a differential $k$-form on $\Omega$. To simplify the notation we assume that $\omega(\bfx) = f(\bfx)dx_1\wedge dx_2\wedge \cdots \wedge dx_k$. Recall that $\mathcal{L}_\bfu\omega = \displaystyle\lim_{t\to 0}\frac{\bfPhi_\bfu(t,\bfx)_*\omega - \omega}{t}$. To simplify the notation, we define $T_t = \bfPhi_\bfu(t)^{-1}$. Then
\begin{align}
\big{(}(\bfPhi_\bfu(t))_*\omega\big{)}(\bfx) = f(T_t(\bfx)) (T_t^*dx_1)  \wedge \cdots  \wedge (T_t^*dx_k) 
\end{align}
It is difficult to have a complete formula for $\bfPhi_\bfu$. We can in turn use its first order Taylor expansion: $\bfPhi_\bfu(t) = \bfx + t\bfu(\bfx) + o(t)$, and $T_t(\bfx) = (\bfPhi_\bfu(t))^{-1} = \bfx - t\bfu(x) + o(t)= (x_1-tu_1(\bfx), x_2-u_2(\bfx), \cdots, x_n-u_n(\bfx)) + o(t)$.  Then direct calculation yields:
\begin{align}
f(T_t(\bfx)) = f(\bfx - t\bfu(\bfx) + o(t)) = f(\bfx) - t\nabla f\cdot\bfu(\bfx) + o(t),
\end{align}
and 
\begin{align}
T_t^*dx_i = d(x_i - tu_i(\bfx) + o(t)) = dx_i - tdu_i +o(t)= dx_i - t\sum_{j=1}^n\frac{\partial u_i}{\partial x_j}dx_j + o(t).
\end{align}
Thus 
\begin{align}
T_t^* \omega &= (f-t\nabla f\cdot \bfu) (dx_1 - tdu_i)\wedge\cdots\wedge (dx_k-tdu_k) + o(t)\\
&= \omega -t\nabla f\cdot\bfu dx_1\wedge\cdots\wedge dx_k - tf\sum_{j=1}^k dx_1\wedge\cdots\wedge dx_{j-1}\wedge du_j\wedge dx_{j+1}\wedge \cdots\wedge dx_k + o(t).
\end{align}
Thus
\begin{align}
\mathcal{L}_{\bfu}\omega = -\nabla f\cdot\bfu dx_1\wedge\cdots\wedge dx_k - f\sum_{j=1}^k dx_1\wedge\cdots\wedge dx_{j-1}\wedge du_j\wedge dx_{j+1}\wedge \cdots\wedge dx_k.
\end{align}

\section{Proof of theorem \ref{thm: existence and uniqueness for tensor fields}}\label{appendix: proof of main theorem}
Recall that $d$ is the exterior derivative operator, $*$ the Hodge star operator on differential forms, and $\delta = (-1)^{n(k-1)+1}*d*$ is the so-called co-differential operator. Let $i:\partial\Omega\rightarrow\Omega$ be the natural embedding of $\partial\Omega$. For convenience but without loss of generality, we set the regularization parameter $a_0 = a_1 = 1$.

First, consider the case when $\Omega$ is a compact oriented Riemannian manifold with boundary. Assume that $\bfu$ is a solution of \eqref{eq: optical flow for tensor fields, compact with boundary, v2}. Take any smooth vector field $\bfh\in T\Omega$ such that $\bfh\big{|}_{\partial\Omega} = 0$. Recall that $\uflat$ and $\hflat$ are differential $1$-forms so that $g(\bfu,\bfv) = \uflat(\bfv)$ and $g(\bfh,\bfv) = \hflat(\bfv)$ for any other vector field $\bfv$. Then we have
\begin{align}
\frac{d}{d \epsilon}\Big{[} \int_{\Omega}2\inp{\theta_1 - \theta_2}{\mathcal{L}_{\bfu+\epsilon \bfh}\theta_2}_\bfx + |d(\uflat+\epsilon \hflat)|_\bfx^2 + |\delta(\uflat+\epsilon \hflat)|_\bfx^2 dV(\bfx) \Big{]}\Big{|}_{\epsilon = 0} = 0
\end{align}
This implies that for any smooth vector field $\bfh\in T\Omega$ such that $\bfh\big{|}_{\partial \Omega} = 0$,
\begin{align}
\int_{\Omega}2\inp{\theta_1 - \theta_2}{\mathcal{L}_\bfh\theta_2}_\bfx + \inp{d\uflat}{d\hflat}_\bfx + \inp{\delta\uflat}{\delta\hflat}_\bfx dV(\bfx) = 0\label{eq: lost function variational form manifold with boundary}
\end{align}
Applying the Green's formula for differential forms (see for instance Eq.(2.1) in \cite{Duff1952HARMONICTO}), we have that
\begin{align}
\int_{\Omega}\inp{d\hflat}{d\uflat}_\bfx dV(\bfx) - \int_{\Omega}\inp{\hflat}{\delta d\uflat}_{\bfx} dV(\bfx) = \int_{\partial\Omega}(i^*\hflat)\wedge *(i^*d\uflat) = 0
\end{align}
since $\hflat\big{|}_{\partial\Omega} = 0$.
Similarly,
\begin{align}
\int_{\Omega}\inp{d\delta\uflat}{\hflat}_\bfx dV(\bfx) - \int_{\Omega}\inp{\delta\uflat}{\delta\hflat}_\bfx dV(\bfx) = \int_{\partial\Omega}(i^*\delta\uflat)\wedge *(i^*\hflat) = 0.
\end{align}

Thus 
\begin{align}
\int_{\Omega}\inp{d\uflat}{d\hflat}_\bfx + \inp{\delta\uflat}{\delta\hflat}dV(\bfx) = \int_{\Omega}\inp{\hflat}{(d\delta+\delta d)\uflat}_\bfx dV(\bfx) = \int_{\Omega}\inp{\hflat}{\Delta\uflat}_\bfx dV(\bfx),\label{eq: second term of variational form}
\end{align}
where $\Delta = d\delta + \delta d$ is the Hodge Laplacian operator.

Let $\alpha = \theta_1-\theta_2, \beta = \theta_2$. By partition of unity,
$\inp{\alpha}{\mathcal{L}_\bfh\beta}$ can be decomposed into a finite sum: 
\begin{align}
\inp{\alpha}{\mathcal{L}_\bfh\beta} = \sum_{i}\inp{\alpha_i}{\mathcal{L}_\bfh\beta_i},
\end{align}
in which $\alpha_i$'s and $\beta_i$'s have the same regularity (hence smooth) as $\alpha$ and $\beta$, and 
are all of the form $\eta_1\otimes\cdots\otimes\eta_l$, with $\eta_j\in V_j$ and $V_j = \Lambda^1(T\Omega)$ or $\Lambda^1(T^*\Omega)$ are smooth sections of the tangent/cotangent bundle respectively. Note that 
\begin{align}
\mathcal{L}_\bfh(\eta_1\otimes\cdots\otimes\eta_l) = \mathcal{L}_\bfh\eta_1\otimes\eta_2\otimes\cdots\otimes\eta_l + \cdots + \eta_1\otimes\cdots\otimes\eta_{l-1}\otimes\mathcal{L}_\bfh\eta_l,
\end{align}
Thus there exists finitely many known smooth functions $m_i(\bfx)$ and known $\xi_i, \eta_i\in \Lambda^1(T\Omega)$ or $\xi_i, \eta_i\in \Lambda^1(T^*\Omega)$, so that 
\begin{align}
\inp{\alpha}{\mathcal{L}_\bfh\beta}_\bfx = \sum_{i}m_i\inp{\xi_i}{\mathcal{L}_\bfh\eta_i}_\bfx.\label{eq: alpha Lh beta = sum mi xi Lh eta}
\end{align}

Since for any $\bfxi,\bfeta\in \Lambda^1(T\Omega)$, $\inp{\bfxi}{\mathcal{L}_\bfh\bfeta} = -\inp{\etaflat}{\mathcal{L}_\bfh\xiflat}$, in which $\xiflat$ (or $\etaflat$) refers to the differential $1$-form such that $\xiflat(\bfv) = \inp{\bfxi}{\bfv}$ (or $\etaflat(\bfv) = \inp{\bfeta}{\bfv}$ resp. ) for any vector field $\bfv\in \Lambda^1(T\Omega)$, without loss of generality we can assume that all the $\xi_i$'s and $\eta_i$'s in Eq. \eqref{eq: alpha Lh beta = sum mi xi Lh eta} are differential $1$-forms. With Cartan's formula $\mathcal{L}_\bfh = d \iota_\bfh + \iota_\bfh d $ and Stokes' formula, we have
\begin{align}
&\int_{\Omega} 2\inp{\theta_1-\theta_2}{\mathcal{L}_{\bfh}\theta_2}_\bfx dV(\bfx) = \int_{\Omega} 2\inp{\alpha}{\mathcal{L}_\bfh\beta}_\bfx dV(\bfx)
= 2\int_{\Omega}\sum_{i}m_i\inp{\xi_i}{\mathcal{L}_\bfh\eta_i}_\bfx dV(\bfx)\nonumber \\
=&2\int_{\Omega}\sum_{i}m_i\mathcal{L}_\bfh\eta_i\wedge *\xi_i = 2\sum_{i}\int_{\Omega}[(d\iota_\bfh + \iota_\bfh d)\eta_i]\wedge *(m_i\xi_i)\nonumber \\
=&2\sum_i\int_{\Omega}d[\iota_\bfh\eta_i\wedge *(m_i\xi_i)] - \iota_\bfh\eta_i\wedge d*(m_i\xi_i) + \iota_\bfh d\eta_i\wedge *(m_i\xi_i)\nonumber \\
=&2\Big{\{}\sum_i\int_{\partial\Omega}\iota_\bfh\eta_i\wedge *(m_i\xi_i) + \int_{\Omega}-\iota_\bfh\eta_i\wedge d*(m_i\xi_i) + \iota_\bfh d\eta_i\wedge *(m_i\xi_i)\Big{\}}\nonumber \\
=& 2\sum_i\int_{\Omega}-\iota_\bfh\eta_i\wedge d*(m_i\xi_i) + \iota_\bfh d\eta_i\wedge *(m_i\xi_i)
\end{align}
For any differential $1$-form $\theta$, denote by $\bfX_\theta\in T\Omega$ the vector field such that $\inp{\eta}{\theta} = \eta(\bfX_\theta)$ for any $\eta\in T^*\Omega$. Here we are not using the standard musical sharp symbol $\theta^{\musSharp{}}$ because in this manuscript $T^{\#}$ denotes how a displacement map $T$ acts on the physical state variables. Then
\begin{align}
&\iota_\bfh d\eta_i\wedge *(m_i\xi_i) = \inp{\iota_\bfh d\eta_i}{m_i\xi_i}dV \nonumber \\
=& (\iota_\bfh d\eta_i)(\bfX_{m_i\xi_i})dV = d\eta_i (\bfh,\bfX_{m_i\xi_i})dV = -d\eta_i(\bfX_{m_i\xi_i},\bfh)dV\nonumber \\
=&-(\iota_{\bfX_{m_i\xi_i}}d\eta_i)(\bfh)dV = -\inp{\iota_{\bfX_{m_i\xi_i}}d\eta_i}{\hflat}dV,
\end{align}
and
\begin{align}
&\iota_\bfh\eta_i\wedge d*(m_i\xi_i) = \inp{\iota_\bfh\eta_i}{*d*m_i\xi_i}dV = \inp{\hflat(\bfX_{\eta_i})}{*d*(m_i\xi_i)}dV\nonumber \\
=& \hflat(\bfX_{\eta_i})*d*(m_i\xi_i)dV = \inp{\hflat}{[*d*(m_i\xi_i)]\eta_i}dV
\end{align}

It thus exists a differential $1$-form $\mu\in \Lambda^1(T^*\Omega)$, which is determined by $\eta_i,m_i,\xi_i$, so that 
\begin{align}
\int_{\Omega}2\inp{\theta_1 - \theta_2}{\mathcal{L}_\bfh\theta_2}_\bfx dV(\bfx) = \int_{\Omega}\inp{\hflat}{\mu}_\bfx dV(\bfx). \label{eq: first term of variational form}
\end{align}
Combining Eq.\eqref{eq: second term of variational form} and \eqref{eq: first term of variational form}, Eq.\eqref{eq: lost function variational form manifold with boundary} is equivalent to
\begin{align}
\int_{\Omega}\inp{\hflat}{\mu}_\bfx + \inp{\hflat}{\Delta \uflat}_\bfx dV(\bfx) = 0.
\end{align}
Thus the optimisation problem \eqref{eq: optical flow for tensor fields, compact with boundary} is equivalent to solving the following equation for $\uflat$:
\begin{align}
\begin{cases}
\Delta \uflat = -\mu\\
\uflat\big{|}_{\partial\Omega} = (\widetilde{i_*\bfv})^{\musFlat{}}\big{|}_{\partial\Omega},\label{eq: laplacian form compact with boundary}
\end{cases}
\end{align}
in which $\widetilde{i_*\bfv}$ is an extension of $\bfv$ from $\partial\Omega$ to $\Omega$ that has the same regularity as $\bfv$ (Theorem 1.3.7 of \cite{Schwarz1995HodgeD}). 
Then the existence, uniqueness, and smoothness of the solution to \eqref{eq: laplacian form compact with boundary} is then guaranteed by theorem (3.4.10) of \cite{Schwarz1995HodgeD}.

For the case, $\Omega$ a compact and oriented Riemannian manifold without boundary, the optimisation problem \eqref{eq: optical flow for tensor fields, compact no boundary} is equivalent to solving the following equation for differential $1$-forms:
\begin{align}
(1-\Delta)\uflat = -\mu,\label{eq: laplacian form compact without boundary}
\end{align}
where $\mu$ is a smooth differential $1$-form that can be derived from the given data. The spectral theory of Laplacian operator on a Riemannian manifold (see for instance theorem (1.30) and (1.31) in  \cite{rosenberg_1997}) states that the space of square integrable differential $k$-forms 
$W^{0,2}(\bigwedge^k T^*\Omega)$ has an orthonormal basis $\{\phi_i: \Delta \phi_i = \lambda_i\phi_i, \displaystyle\int_{\Omega}\inp{\phi_i}{\phi_i}_\bfx dV(\bfx) = 1,  0\leq i <\infty\}$, and that all eigen-forms $\phi_i$ are smooth on $\Omega$. Note that $\lambda_i \leq 0$ $\forall i$. Thus we can assume the decomposition:
$\mu = \sum_{i\geq 0}a_i\phi_i$.
Then $\uflat = \displaystyle\sum_{i}\frac{a_i}{1-\lambda_i}\phi_i$ is a solution to Eq.\eqref{eq: laplacian form compact without boundary}. The identity can be verified directly. In addition we need to show that $\uflat$ is smooth, i.e $\uflat$ has finite $H^s$ norm for any integer $s \geq  0 $. This can be proved by induction based on Theorem 2.46 in \cite{rosenberg_1997}. For $s=0$, the finiteness of $\|\uflat\|_0$ can be directly verified:
\begin{align}
\|\uflat\|_0^2 = \sum_{i}(\frac{a_i}{1-\lambda_i})^2\leq \sum_{i}|a_i|^2 = \|\mu\|_0^2.
\end{align}
G\r{a}rding's inequality (theorem 2.46 in \cite{rosenberg_1997}) shows that, on a compact manifold without boundary, the $H^1$ norm of any differential form $\theta$ is equivalent to the Dirichlet integral: $\|d\theta\|_0^2 + \|\delta \theta\|_0^2 + \|\theta\|_0^2$. It is not hard to show that the infinite sum $\sum_{i=0}^{\infty}a_i\psi_i$ converges under the norm defined by Dirichlet integral. Thus $\uflat\in H^1$ and  $d\uflat, \delta\uflat\in H^0$, hence $D\uflat: = (d+\delta)\uflat\in H^0$. Then  applying the theorem 2.46 in \cite{rosenberg_1997}, for any vector field $\bfv$ we have that
\begin{align}
\|\vflat\|_2\lesssim \|\vflat\|_1 + \|D\vflat\|_1\lesssim \|\vflat\|_1 + \|D\vflat\|_0 + \|dD\vflat\|_0 + \|\delta D\vflat\|_0 =  \|\vflat\|_1 + \|D\vflat\|_0 + \|d\delta \vflat\|_0 + \|\delta d\vflat\|_0
\end{align}

Now take $\vflat_n = \displaystyle\sum_{i=1}^n\frac{a_i}{1-\lambda_i}\psi_i$, then
\begin{align}
&\|d\delta\vflat_n\|_0^2 + \|\delta d\vflat_n\|_0^2= \int_{\Omega}\inp{d\delta\vflat_n}{d\delta\vflat_n}_{\bfx} + \inp{\delta d\vflat_n}{\delta d\vflat_n}_{\bfx}dV(\bfx)\\
=& \int_{\Omega}\inp{\Delta^2\vflat_n}{\vflat_n}_{\bfx}dV(\bfx) = \sum_{i=1}^n\frac{\lambda_i^2}{(1-\lambda_i)^2}a_i^2\leq \sum_{i}a_i^2 = \|\mu\|_0^2.
\end{align}
Similar arguments show that $\{\vflat_n\}$ is a Cauchy sequence in $H^2$ thus  $\uflat\in H^2$. Thus $\Delta \uflat$ is well-defined and $(I-\Delta)\uflat = -\mu$ indeed holds. Now we continue to show that $\uflat$ is smooth. Suppose that we have shown that $\uflat\in H^{s}$ for some integer $s\geq 2$. By applying Theorem 2.46 in \cite{rosenberg_1997} we have:
\begin{align}
&\|\uflat\|_{s+1}\lesssim \|\uflat\|_s + \|D\uflat\|_{s}\lesssim \|\uflat\|_s + \|D\uflat\|_{s-1} + \|DD\uflat\|_{s-1} = \|\uflat\|_s + \|D\uflat\|_{s-1} + \|\Delta \uflat\|_{s-1}\\
=& \|\uflat\|_s + \|D\uflat\|_{s-1} + \|\uflat + \mu\|_{s-1} \leq \|\uflat\|_s + \|D\uflat\|_{s-1} + \|\uflat \|_{s-1} + \|\mu\|_{s-1} <\infty.
\end{align}

Thus $\uflat\in H^{s+1}$. Thus $\uflat\in H^s$ for all $s \geq 0$.

\section{Proof of theorem \ref{thm: OF of tensor fields}}
First assume that the manifold has $C^1-$boundary. Then the boundary condition $\bfv$ can be extended to $\bfu_0\in H^1(T\Omega)$ such that $\bfu_0\big{|}_{\partial\Omega} = \bfv$. Thus the optimisation problem \eqref{eq: OF of tensor fields, with boundary} is equivalent to 
\begin{align}
\begin{cases}
\bfu(t) = \displaystyle\argmin_{\bfu\in H^1(T\Omega)}\int_{\Omega}|\theta_t - \mathcal{L}_{\bfu_0}\theta - \mathcal{L}_{\bfu}\theta|_\bfx^2 + |d\uflat + d\uflat_0|_\bfx^2 + |\delta \uflat + \delta\uflat_0|_\bfx^2 \\
\bfu(t)\big{|}_{\partial\Omega} = 0.
\end{cases}\label{eq: OF of tensor fields, with boundary, v2}
\end{align}
The first term in the above functional can be rewritten as 
\begin{align}
\inp{\mathcal{L}_\bfu\theta}{\mathcal{L}_\bfu\theta}_\bfx - 2\inp{\mathcal{L}_\bfu\theta}{\theta_t - \mathcal{L}_{\bfu_0}\theta}_\bfx + |\theta_t - \mathcal{L}_{\bfu_0}\theta|_\bfx^2
\end{align}
Since $\theta\in H^1$ and $\bfu_0\in H^1$, $\mathcal{L}_{\bfu_0}\theta\in L^2$. Note that for the space of vector fields vanishing on $\partial\Omega$, the $L^2$ norm of $\bfu$ is bounded by $|d\uflat|^2 + |\delta\uflat|^2$ up to a constant depending on the domain only.  For $\bfu,\bfw$ of finite $H^1$ norm and vanishing on the boundary, let
\begin{align}
B(\bfu,\bfw) &= \int_{\Omega}\inp{\mathcal{L}_\bfu\theta}{\mathcal{L}_\bfw\theta}_x + \inp{d\uflat}{d\wflat}_x + \inp{\delta\uflat}{\delta\wflat}_\bfx dV(\bfx)\\
a(\bfu) &= \int_{\Omega}-\inp{\mathcal{L}_\bfu\theta}{\theta_t-\mathcal{L}_{\bfu_0}\theta}_\bfx + \inp{d\uflat}{d\uflat_0}_{\bfx} + \inp{\delta\uflat}{\delta\uflat_0}_\bfx dV(\bfx)
\end{align}
Obviously $B$ is symmetric. $B$ is further coercive due to Poincare lemma and Gaffney's inequality (See equation (0.9) of \cite{Schwarz1995HodgeD} for the case when $\Omega$ is compact with boundary). We will show that $B$ is bounded and $a$ is continuous with respect to the $H^1$ norm. Then $B$ gives an equivalent norm as the common $H^1$ norm, denoted by $\|\cdot\|_B$.  Then by Riesz representation theorem, there exists $\bfu_1\in H^1(T\Omega)$ which vanishes on $\partial\Omega$, such that $a(\bfu) = \inp{\bfu}{\bfu_1}_B = B(\bfu,\bfu_1)$ for any $\bfu\in H^1(T\Omega)$ that vanishes on $\partial\Omega$. Then
\begin{align}
B(\bfu,\bfu) + a(\bfu) = \inp{\bfu}{\bfu}_B - 2\inp{\bfu}{\bfu_1}_B = \inp{\bfu-\bfu_1}{\bfu-\bfu_1}_B - \|\bfu_1\|_B^2\geq -\|\bfu_1\|_B^2,
\end{align}
implying that $\bfu = \bfu_1$ is the unique solution. Next we show that $B$ and $a$ are continuous forms.

For the continuity of $a$, following a similar argument as those in appendix \ref{appendix: proof of main theorem}, there exists finitely many $\eta_i\in H^1(T^*\Omega)$ with bounded $|d\eta_i|_\bfx$ and $\xi_i\in L^2(T^*\Omega)$ depending only on the given data, such that
\begin{align}
a(\bfu) = \sum_{i}\int_{\Omega} \inp{\mathcal{L}_\bfu\eta_i}{\xi_i}_\bfx dV(\bfx)
\end{align}
We can further assume that $\eta_i,\xi_i$ are compactly supported inside a local coordinate $(\Omega_i, \bfx)$. Since $|d\eta_i|_\bfx$ is bounded, 
\begin{align}
\|\mathcal{L}_\bfu\eta_i\|_{L^2} = \|d\iota_\bfu\eta_i + \iota_\bfu d\eta_i\|_{L^2} = \|d(\eta_i(\bfu)) + \iota_\bfu d\eta_i\|_{L^2}\lesssim \|\bfu\|_{H^1}.
\end{align}
Thus 
\begin{align}
|a(\bfu)| \lesssim \sum_{i}\|\mathcal{L}_\bfu\eta_i\|_{L^2}\|\xi_i\|_{L^2}\lesssim \|\bfu\|_{H^1}.
\end{align}

In order to show that $B$ is bounded, following the same argument in appendix \ref{appendix: proof of main theorem}, there exists finitely many differential $1$-forms $\eta_i$ in $H^1$ with bounded $|d\eta_i|_\bfx$, and bounded functions $m_i$, such that 
\begin{align}
\inp{\mathcal{L}_\bfu\theta}{\mathcal{L}_\bfu\theta}_\bfx = \sum_{i}m_i\inp{\mathcal{L}_\bfu\eta_i}{\mathcal{L}_\bfu\eta_i}_\bfx.
\end{align}
Therefore  $\|\mathcal{L}_\bfu\theta\|^2_{L^2}\lesssim \|\bfu\|_{H^1}^2$, meaning that $B$ is bounded. The proof for the case when $\Omega$ is a compact oriented Riemannian manifold without boundary is similar, thus omitted.

\section{On the linearity of $\mathcal{L}_{\bfu(t)}\theta_2$ with respect to $\theta_1 - \theta_2$}\label{appendix: the linearity of L_u theta on y-xb}
In this appendix we assume that quantities $\theta_1$, $\theta_2$, and $\bfu(t)$ are the same as defined in Eq.\eqref{eq: optical flow for tensor fields, compact no boundary, v2} or \eqref{eq: optical flow for tensor fields, compact with boundary, v2}, and the quantity $\mu$ is the one defined in the proof of theorem 2.1 in appendix \ref{appendix: proof of main theorem}. The goal here is to show that $\mathcal{L}_{\bfu(t)}\theta_2$ is linear with respect to $\theta_1 - \theta_2$. In the proof of theorem 2.1 (see appendix \ref{appendix: proof of main theorem}), it has been proved that there exists a differential  $1$-form $\mu$, so that $\displaystyle\int_{\Omega}\inp{\theta_1-\theta_2}{\mathcal{L}_\bfh\theta_2}_\bfx dV(\bfx)= \int_{\Omega}\inp{\hflat}{\mu}_\bfx dV(\bfx)$ for any $\bfu$ which vanishes on the boundary of $\Omega$. The construction of $\mu$ depends on the partition of unity and the operators $d$, $*$, $\iota$, and the musical sharp and musical flat operators. Thus we can write $\mu = l(\theta_2, \theta_1 - \theta_2)$, for some operator $l(\cdot, \cdot)$ which does not depend on $\bfh$. Denote by $\alpha = \theta_2$ and $\beta = \theta_1 - \theta_2$, we want to show that $l(\alpha, \beta)$ is bilinear. This is obvious because every equation is bilinear in $\alpha$ and $\beta$ in the construction of $\mu$ in appendix \ref{appendix: proof of main theorem}. Thus $\uflat = \Delta^{-1}\mu$ or $(I-\Delta)^{-1}\mu$ is linear in $\theta_1 - \theta_2$. And because $\mathcal{L}_{\bfu}\theta$ is linear in $\bfu$, it is also linear in $\theta_1 - \theta_2$.

\section{Details of the numerical experiments}
\subsection{Model and domain}
The data assimilation experiment is conducted using the thermal shallow water equation \cite{Warneford2013TheQT}. This model consists of three state variables: $h-$the water height, $\bfv = (v_1,v_2)-$the velocity field, and $\Theta-$the buoyancy (or density contrast):
\begin{align}
\frac{\partial h}{\partial t} + \nabla\cdot(h \bfv) &= 0,\\
\frac{\partial \Theta}{\partial t} + (\bfv\cdot\nabla)\Theta &= -\kappa(h\Theta - h_0\Theta_0),\\
\frac{\partial \bfv}{\partial t} + (\bfv\cdot\nabla)\bfv + f\hat{\bf z}\times \bfv &= -\nabla(h\Theta) + \frac{1}{2}h\nabla\Theta.
\end{align}
Both $\Theta$ and $h$ are assumed strictly positive at each point. 

The quantity $h$ follows a dynamics similar to the state variable $S$ in example \ref{example: transport eqn, n-form}, and transportation terms for $\Theta$ follow those in example \ref{example: transport eqn, 0-form}. A natural choice is then $\theta_h = h dx_1\wedge dx_2$, a differential 2-form, and $\theta_\Theta = \Theta$ a differential 0-form. The test case is a double-vortex case. Hence $\theta_\bfv = v_1dx_1 + v_2dx_2$ is chosen to ensure vorticity conservation during the morphing process. These choices of differential forms differ from those presented in \cite{Zhen2023PhysicallyCC}. Less constrained by the underlying dynamics, \cite{Zhen2023PhysicallyCC} discussed the derivation of a perturbation scheme able to conserve particular quantities. Here, the choice for the tensor fields should obey to the prescribed dynamics of the system to maximally maintain the dynamical balance during the morphing process.

The data assimilation process is solely conducted for one time step. It is assumed that both the vorticity field $\omega = \frac{\partial v^2}{\partial x^1} - \frac{\partial v^1}{\partial x^2}$ and the $h$ field are fully observed. Since $\theta_\bfv = v_1dx_1+v_2dx_2$, naturally $\omega$ is associated to a differential 2-form $\theta_\omega = \omega dx_1\wedge dx_2 = d\theta_\bfv$.

The domain is 2-dimensional doubly periodic: $\Omega = [0,5000\text{km}]_{\text{per}}\times [0,5000\text{km}]_{\text{per}}$, which is a compact Riemannian manifold without boundary. In this case, a numerical solution of \eqref{eq: optical flow for tensor fields, compact no boundary, v2}, or equivalently the numerical solution of Eq.\eqref{eq: laplacian form compact without boundary}, with $a = 1$, can be derived in the Fourier space. In fact, the vector field $u$ is separately calculated for $h$ observations and for $\omega$ observations. From the two observables, the explicit expressions of $u$ are:
\begin{align}
\bfu_{\omega} &= (u_{\omega,1}, u_{\omega,2}) = 2(I-\Delta)^{-1} \big{[}\omega_2\nabla(\omega_1-\omega_2)\big{]},\label{eq: u-omega} \\
\bfu_{h} &= (u_{h,1}, u_{h,2}) = 2(I-\Delta)^{-1} \big{[}h_2\nabla(h_1-h_2)\big{]}.\label{eq: u-h}
\end{align}
The final $\bfu$ for each iterative step is then chosen to be $\frac{1}{2}(\frac{\bfu_\omega}{\|\bfu_{\omega}\|_1}+\frac{\bfu_h}{\|\bfu_h\|_1})$, in which 
\begin{align}
\|\bfu\|_1^2 = \int_{\Omega} \inp{\uflat}{\uflat}_\bfx + \inp{d\uflat}{d\uflat}_\bfx + \inp{\delta\uflat}{\delta\uflat}_\bfx dV(\bfx) =\displaystyle\int_{\Omega}|\bfu|^2 + (\frac{\partial u_2}{\partial x_1} - \frac{\partial u_1}{\partial x_2})^2 + (\frac{\partial u_{1}}{\partial x_1} + \frac{\partial u_2}{\partial x_2})^2 dx_1dx_2.
\end{align}

\subsection{Numerical methods and experimental parameters}
The units in the code are set to be km and 100s. The initial condition and model parameters are taken from the numerical experiment in subsection (5.3) of \cite{Eldred2019AQD}. The 3-step Adams-Bashforth method (see for instance chapter 3.1 of \cite{Hairer1993} ) is used in model integration. Additionally, the one dimensional Hou-Li spectral filter \cite{Hou2007ComputingNS} $\exp\{-36[(k_x/k_{max})^a + (k_y/k_{max})^a ]\}$, with $a = 12$,  is applied to the Fourier modes of the state vector at the end of each model integration step. The truth is generated by running the model forward for 2750 time units. To generate the ensemble members, the center of the initial vortex, $(ox,oy)$, is perturbed:
\begin{align}
ox \sim \mathcal{N}(0.1,0.01),\hspace{2mm}
oy \sim \mathcal{N}(0.1,0.01).
\end{align}
The ensemble members are then generated by running the model forward starting from perturbed initial condition for 2000 time units. For the explicit meaning of $ox$ and $oy$, please refer to section 5.3 of \cite{Eldred2019AQD}. The ensemble size $N_e = 20$. 
Both the ensemble members and the truth are generated using a $256\times 256$ grid. But before starting the morphing process or data assimilation, the ensemble members and the observations are both projected to a coarse-grid (64 $\times$ 64) and then interpolated back to the original grid ($256\times 256$).

The morphing process is implemented with the 5-step Adams-Bashforth method and the spectral Hou-Li filter with $a = 36$. For clarity, the pseudo-code of the complete morphing process is shown in Algorithm 2.%\ref{alg: morphing process for thermal shallow water exp}. 
We choose $\epsilon = 0.000033$ and $N = 10000$.

\begin{figure}[h]
	\centering
	\includegraphics[width=0.9\linewidth]{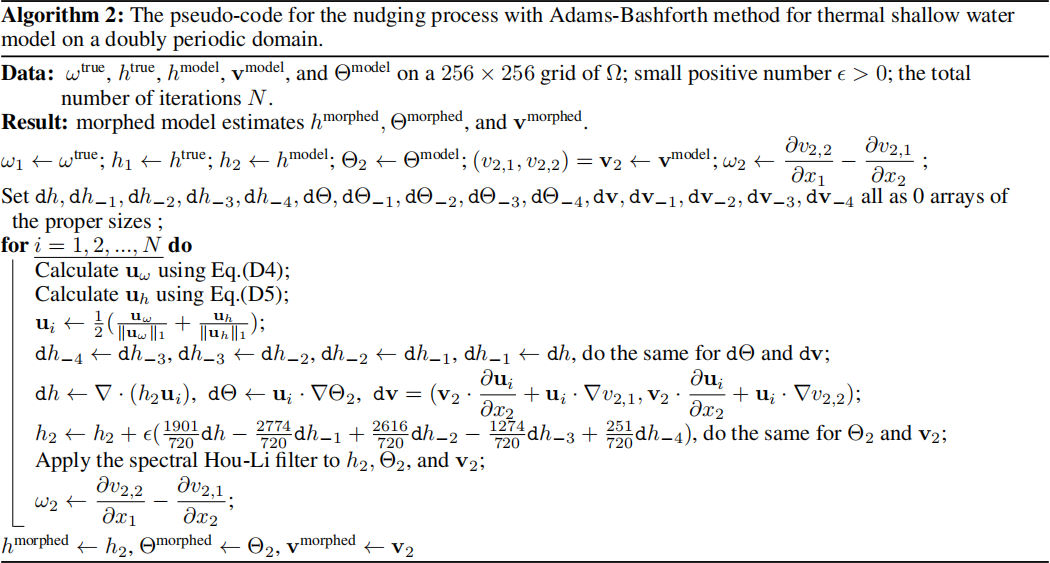}
	\label{alg: morphing process for thermal shallow water exp}
\end{figure}

The observation is exactly the same as the interpolated value of the truth. There is no error in the observation. However, the matrix $R$ in ensemble Kalman filter is set to be a diagonal matrix with diagonal elements equal to
\begin{align}
R_{\omega} = \frac{0.01}{64^2}\sum_{i=1}^{64^2}(\omega^{\text{obs}}_i)^2, \text{\hspace{1mm} or \hspace{1mm}} R_{h} = \frac{0.01}{64^2}\sum_{i=1}^{64^2}(h^{\text{obs}}_i)^2,
\end{align}
where $\omega_i^{\text{obs}}$ and $h_i^{\text{obs}}$ refer to the $\omega$ value and $h$ value at the $i-$th grid-point in the $64\times 64$ grid. Data assimilation is conducted on the 64$\times$64 grid.

%% Please add \clearpage between each table and/or figure. Further guidelines on figures and tables can be found below.

%\disclaimer{TEXT} %% optional section

\section*{Acknowledgements}
	The authors would like to express their gratitude to Long Li at the department of mathematics of ShanghaiTech University and Wei Pan at department of mathematics of Imperial College London for their insightful and helpful discussion, the comments from both anonymous reviewers and Etienne M\'{e}min for his introduction of relevant works. YZ's research is supported by National Natural Science Foundation of China (NSFC, Grant No. 42350003). VR's research is supported by the company SCALIAN DS and by France Relance through the MORAANE project. BC's research is supported by ERC EU SYNERGY Project No. 856408-STUOD, and the the support of the ANR Melody project.

%\printbibliography[
%heading=bibintoc,
%title={References}
%] 

\end{document}